\documentclass[11pt,reqno]{amsart}

\usepackage{graphics,enumitem,epsfig,textcomp,bm}
\usepackage{amsfonts,amsmath,amssymb,euscript,color,mathrsfs}

\numberwithin{equation}{section}
\newtheorem{defin}{Definition}%[section]
\newtheorem{theorem}{Theorem}%[section]

\newtheorem{lemma}{Lemma}%[section]
\def\begproof{\noindent{\bf Proof: }}
\def\endproof{\par\rightline{\vrule height5pt width5pt depth0pt}\medskip}
\def\d{\,\mathrm{d}}
\def\eps{\varepsilon}
\def\phi{\varphi}
\def\theta{\vartheta}
\def\N{\mathbb{N}}
\def\R{\mathbb{R}}
\def\C{\hbox{\rlap{\kern.24em\raise.1ex\hbox
      {\vrule height1.3ex width.9pt}}C}}
\def\P{\hbox{\rlap{I}\kern.16em P}}
\def\Q{\hbox{\rlap{\kern.24em\raise.1ex\hbox
      {\vrule height1.3ex width.9pt}}Q}}
\def\M{\hbox{\rlap{I}\kern.16em\rlap{I}M}}
\def\Z{\hbox{\rlap{Z}\kern.20em Z}}
\def\({\begin{eqnarray}}
\def\){\end{eqnarray}}
\def\[{\begin{eqnarray*}}
\def\]{\end{eqnarray*}}
\def\part#1#2{\frac{\partial #1}{\partial #2}}
\def\partk#1#2#3{\frac{\partial^#3 #1}{\partial #2^#3}}

\def\grad{\nabla}
\def\Norm#1{\left\| #1 \right\|}
\def\pmb#1{\setbox0=\hbox{$#1$}
  \kern-.025em\copy0\kern-\wd0
  \kern-.05em\copy0\kern-\wd0
  \kern-.025em\raise.0433em\box0 }
\def\bar{\overline}

\def\mmnt#1{\left\langle #1 \right\rangle}

\def\tot#1#2{\frac{\d #1}{\d #2}} 
\def\laplace{\Delta}
\def\d{\,\mathrm{d}}
\def\N{\mathbb{N}}
\def\R{\mathbb{R}}

\def\eps{\varepsilon}

\def\AA{\mathcal{A}}

\def\bust{\bar u^*}
\def\xe{\xi_\eps}
\def\Ae{A_\eps}
\def\Fe{F_\eps}
\def\De{D_\eps}
\def\we{{\bf w}_\eps}
\def\uei{u^\eps_\infty}

\def\bA{{\bf A}}
\def\bB{{\bf B}}
\def\bx{{\bf x}}
\def\by{{\bf y}}
\def\balpha{{\bm{\alpha}}}

\def\aA{\balpha\cdot\bA}
\def\AA{\mathcal{B}}

\def\H{\mathcal{H}}
\def\L{\mathcal{L}}

\def\comment#1{\textcolor{blue}{\bf [#1]}}

\newcommand{\md}{d\kern-0.035cm\char39\kern-0.03cm}

\title[Well posedness and approximation for quantitative traits]{Well posedness and Maximum Entropy Approximation for the Dynamics of Quantitative Traits}

\author[Bo\v dov\' a]{Katar\'\i na Bo{\v d}ov\' a}
\address[Katar\'\i na Bo{\md}ov\' a]{\newline Institute of Science and Technology Austria (IST Austria), Klosterneuburg A-3400, Austria}
\email{kbodova@ist.ac.at}

\author[Haskovec]{Jan Haskovec}
\address[Jan Haskovec]{\newline Computer, Electrical and Mathematical Sciences \& Engineering
    \newline King Abdullah University of Science and Technology, 23955 Thuwal, KSA}
\email{jan.haskovec@kaust.edu.sa}

\author[Markowich]{Peter Markowich}
\address[Peter Markowich]{\newline Computer, Electrical and Mathematical Sciences \& Engineering
    \newline King Abdullah University of Science and Technology, 23955 Thuwal, KSA}
\email{peter.markowich@kaust.edu.sa}

\begin{document}
%%%%%%%%%%%%%%%%
\allowdisplaybreaks

\date{\today}

\subjclass[]{}

\begin{abstract} We study the Fokker Planck equation derived in the large system limit
of the Markovian process describing the dynamics of quantitative traits.
The Fokker-Planck equation is posed on a bounded domain and its
transport and diffusion coefficients vanish on the domain's boundary.
We first argue that, despite this degeneracy, the standard no-flux boundary condition
is valid. We derive the weak formulation of the problem
and prove the existence and uniqueness of its solutions
by constructing the corresponding contraction semigroup
on a suitable function space.
Then, we prove that for the parameter regime with high enough mutation rate
the problem exhibits a positive spectral gap, which implies exponential
convergence to equilibrium.

Next, we provide a simple derivation of the so-called Dynamic Maximum Entropy (DynMaxEnt)
method for approximation of moments of the Fokker-Planck solution,
which can be interpreted as a nonlinear Galerkin approximation.
The limited applicability of the DynMaxEnt method
inspires us to introduce its modified version that is valid for the whole range of admissible parameters.
Finally, we present several numerical experiments to demonstrate the performance
of both the original and modified DynMaxEnt methods. We observe that in the parameter regimes
where both methods are valid, the modified one exhibits slightly better approximation properties
compared to the original one.
\end{abstract}

\maketitle \centerline{\date}

\tableofcontents

%%%%%%%%%%%%%%%%
%   Introduction
%%%%%%%%%%%%%%%%
\section{Introduction}
The dynamics of allele frequencies $\bx = (x_1, \dots, x_L)$, 
where $L$ is the number of loci that contribute to the trait,
can be described by a diffusion process using a deterministic forward Kolmogorov equation.
The evolution of the joint probability density $u = u(t,\bx)$ of allele frequencies
for a population of $N$ diploid individuals satisfies the linear Fokker-Planck equation
\( \label{eq0}
   \part{u}{t} = - \frac12 \sum_{i=1}^L \part{}{x_i} \left( \xi_i \part{(\aA)}{x_i} u \right)
      + \frac1{4N} \sum_{i=1}^L \partk{}{x_i}{2} (\xi_i u),
\)
on $\Omega_\bx := (0,1)^L$,
where we denoted $\xi_i := \xi(x_i) = x_i(1-x_i)$ for $i=1,\dots, L$.
The diffusion term captures the stochasticity of the allele frequencies
arising from random sampling. %, i.e., random drift.
Here we assume that linkage disequilibria are negligible,
otherwise this term would be of cross-diffusion type,
reflecting correlations between loci \cite{BTB}.
The drift term captures deterministic effects on allele frequencies
that are described by a vector of coefficients $\balpha$
and a vector of complementary quantities $\bA$.
We consider directional selection and dominance with symmetrical mutation,
which, using the notation of \cite{BTB}, corresponds to the choice
\[
   \bA = (\xi'_1, \dots, \xi'_L, \xi_1,\dots,\xi_L, \ln\xi_1,\dots,\ln\xi_L)
\]
and
\(  \label{alphaA}
   \aA = - \beta\sum_{i=1}^L \gamma_i \xi_i' + 2h \sum_{i=1}^L \eta_i \xi_i + 2\mu \sum_{i=1}^L \ln \xi_i,
\)
where the nondimensional parameters $\beta, h, \gamma_i, \eta_i \in\R$ represent the effects of loci on the traits,
$\mu>0$ is the mutation rate, and $\xi_i' := \xi'(x_i) = 1-2x_i$.
For notational simplicity and without loss of generality, we set $\beta = h = 1$ in the sequel, so that
\[
   \balpha = (-\gamma_1, \dots, -\gamma_L, 2\eta_1,\dots,2\eta_L, 2\mu,\dots,2\mu)\in\R^{3L} \,.
\]

This drift-diffusion process \eqref{eq0} is known to be an accurate
continuous-time approximation to a wide range of specific population
genetics models \cite{Kimura55a, Kimura55b, Ewens12, dV-Barton12}.
%Moreover, it corresponds directly to the coalescent process that describes the ancestry of samples taken from the population \cite{Wakeley08}.
In order to represent the population
in terms of allele frequencies, we must assume that linkage
disequilibria are negligible, which will be accurate if recombination
is sufficiently fast. For simplicity, we also assume two
alleles per locus.

The main difficulty for analysis of the Fokker-Planck equation \eqref{eq0}
is the degeneracy of the diffusion coefficients $\xi_i = x_i(1-x_i)$ at the boundary
of $\Omega_\bx$. Consequently, the task of prescribing boundary conditions
that lead to a well-posed problem is far from obvious; see also \cite{Chalub1, Chalub2}
for related issues in population genetics problems.
As noted above, we aim at interpreting the solution $u$ as a time-dependent probability density,
which calls for a no-flux boundary condition.
In Section \ref{sec:BC} we argue that the standard no-flux boundary condition is indeed
appropriate for \eqref{eq0}. In Section \ref{sec:WP} we derive the weak formulation of \eqref{eq0}
subject to the no-flux boundary condition and prove the existence and uniqueness of its solutions
by constructing the corresponding contraction semigroup.
Then, in Section \ref{sec:gap} we prove that for the parameter regime with high enough mutation rate
the problem exhibits a positive spectral gap, which implies exponential
convergence to equilibrium.

In typical applications in quantitative genetics the solution of the
Fokker-Planck equation \eqref{eq0} is not the main object of interest.
One is rather interested in the evolution of its certain moments
that correspond to the macroscopic dynamics of observable quantitative traits.
Therefore, Section \ref{sec:DNE} is devoted to the study of the so-called Dynamic Maximum Entropy (DynMaxEnt)
method for approximation of moments of the Fokker-Planck solution.
We first show in Section \ref{subsec:entropy} that a related constrained entropy maximization is equivalent
to a moment-matching problem, which we solve in a simple case.
Then, in Section \ref{subsec:DNE} we provide a simple and straightforward derivation
of the DynMaxEnt method by adopting a quasi-stationary approximation,
which results in a nonlinear system of ordinary differential equations.
It can be interpreted as a nonlinear Galerkin approximation of the Fokker-Planck
equation \eqref{eq0}.
However, this "original" DynMaxEnt method cannot be applied in the regime of small mutations, i.e., when $4N\mu \leq 1$.
This inspires us to introduce a modified version, which is valid for the whole range of admissible parameters,
Section \ref{subsec:simple}.
Finally, in Section \ref{sec:numerics} we present several numerical experiments to demonstrate the performance
of both the original and modified DynMaxEnt methods. We observe that in the parameter regimes
where both methods are valid, the modified one exhibits slightly better approximation properties
compared to the original one.

The surprisingly good approximation properties of the DynMaxEnt method,
as documented by the numerical results in \cite{BTB} and Section \ref{sec:numerics}
of this paper, suggest that the infinitely-dimensional dynamics of the Fokker-Planck equation \eqref{eq0}
can be well approximated by suitable finitely-dimensional dynamical systems.
This is reminiscent of the recent series of works of E. Titi and collaborators \cite{MT, FLT, ATKZ, MTT, ALT}
where a data assimilation (downscaling) approach to fluid flow problems is developed,
inspired by ideas applied for designing finite-parameters feedback control for dissipative systems.
The goal of a data assimilation algorithm is to obtain (numerical) approximation
of a solution of an infinitely-dimensional dynamical system corresponding to given measurements
of a finite number of observables. In particular, in \cite{MT}, it has been shown that solutions 
of the two-dimensional Navier-Stokes equations can be well reconstructed from
a relatively low number of low Fourier modes or local averages over finite volume elements.
In \cite{FLT}, continuous data assimilation (CPA) algorithm was proposed and analyzed for a two-dimensional B\'enard convection problem,
where the observables were incorporated as a feedback (nudging) term in the evolution equation of the horizontal velocity.
In \cite{ATKZ} CPA was applied for downscaling a coarse
resolution configuration of the 2D B\'enard convection equations into a finer grid,
while in \cite{MTT} the CPA method is studied for a three-dimensional
Brinkman-Forchheimer-extended Darcy model of porous media,
and in \cite{ALT} for the three-dimensional Navier-Stokes--$\alpha$ model.
Finally, in \cite{GOT} numerical performance of the CPA algorithm in the context
of the two-dimensional incompressible Navier--Stokes equations was studied.
It was shown that the numerical method is computationally efficient
and performs far better than the analytical estimates suggest.
This is similar to our numerical observations showing
very good approximation properties of the DynMaxEnt method
applied to the Fokker-Planck equation \eqref{eq0}.

%%%%%%%%%%%%%%%%%%%%%%%%%%%%%%%%%
\section{Boundary conditions for the stationary problem}\label{sec:BC}
The stationary solution of the Fokker-Planck equation \eqref{eq0} is of the form
\(   \label{stat0}
   u_\balpha = \frac{1}{\mathbb Z_\balpha} \frac{\exp(2N\aA)}{\prod_{i=1}^L \xi_i},
\)
where $\mathbb Z$ is a normalization constant (partition function).
We aim at interpreting the solution $u_\balpha$ as a probability density,
therefore, we set
\(   \label{Z_alpha}
   \mathbb Z_\balpha := \int_{\Omega_\bx} \frac{\exp(2N\aA)}{\prod_{i=1}^L \xi_i} \d x.
\)
Observe that the above integral is finite for $\mu>0$, which we assumed.
Let us rewrite \eqref{eq0} in the form
\( \label{eq1}
   \part{u}{t} = \grad_\bx\cdot \left( D u_\balpha \grad_\bx \left(\frac{u}{u_\balpha}\right)\right)
\)
with $u_\balpha$ defined in \eqref{stat0}, and
the diagonal diffusion matrix $D = D(\bx)$, $D_{ij} = \frac{1}{4N} \xi_i  \delta_{ij}$ for $i,j=1,\dots,L$.

To provide an insight into the problem of prescribing valid boundary conditions
for \eqref{eq1}, we consider the related stationary problem in the spatially one-dimensional setting,
\( \label{stationary1D}
   \partial_x \left( D u_\balpha \partial_x \left(\frac{u}{u_\balpha}\right)\right) = f
\)
for $x\in (0,1)$, where $f\in L^1(0,1)$ is a prescribed function with $\int_0^1 f(s) \d s = 0$,
$\xi = \xi(x) = x(1-x)$, $D = \frac{1}{4N}\xi$ and
\[
   u_\balpha = \mathbb Z_\balpha^{-1} \xi^{-1}\exp (2N\aA) =  \mathbb Z_\balpha^{-1} \xi^{4N\mu -1}\exp (2N\gamma\xi' + 4N\eta\xi),
\]
with $\mathbb Z_\balpha$ defined in \eqref{Z_alpha}.
We recall that $\mathbb Z_\balpha$ is finite and $u_\balpha$ is integrable for the relevant range of parameters.
Moreover, note that the product $D u_\balpha$ behaves like $\xi^{4N\mu}$ close to $x = 0$ and $x=1$,
so that it vanishes at the boundary and leads to a degeneracy in the formal no-flux boundary condition
\(  \label{noFlux1D}
   D u_\balpha \partial_x \left(\frac{u}{u_\balpha}\right) = 0 \quad\mbox{for } x\in\{0,1\}. 
\)
To avoid possible difficulties due to this degeneracy,
we integrate \eqref{stationary1D} for a fixed $x\in(0,1)$ on the interval $(1/2,x)$,
\[
   D u_\balpha \partial_x \left(\frac{u}{u_\balpha}\right) = \int_{1/2}^x f(s) \d s + C_1,
\]
where $C_1$ is an integration constant.
We see that imposing the formal no-flux boundary condition \eqref{noFlux1D} at, say, $x=0$
is equivalent to setting $C_1$ to the particular value
\[
   C_1 = \int_0^{1/2} f(s) \d s.
\]
The assumption $\int_0^1 f(s) \d s = 0$ then implies that \eqref{noFlux1D} is verified at $x=1$.
Integrating once again yields
\(  \label{u_formula}
   u = C_2 u_\balpha + C_1 u_\balpha \int_{1/2}^x \frac{\d s}{D(s) u_\balpha (s)}
      + u_\balpha \int_{1/2}^x \frac{F(s)\d s}{D(s) u_\balpha(s)},
\)
with $F(s) := \int_{1/2}^s f(r) \d r$. Observe that
\[
   \int_{1/2}^x  \frac{\d s}{D(s) u_\balpha (s)} \approx \xi^{-4N\mu + 1} \quad\mbox{close to } x\in\{0,1\},
\]
so that the second term in \eqref{u_formula} is bounded on $[0,1]$ and thus integrable.
Due to the boundedness of $F(s)$, the same holds also for the third term in \eqref{u_formula}.
Consequently, the solution $u$ constructed in \eqref{u_formula} is integrable on $(0,1)$.

We conclude that, for the aforementioned range of parameter values,
the Fokker-Planck equation \eqref{eq1} has to be supplemented with the standard
no-flux boundary condition \eqref{noFlux1D} regardless of the degeneracy of $D u_\balpha$ at the boundary.
Although the above argument only applies to the spatially one-dimensional setting,
it provides a strong heuristic hint that the conclusion also holds in the multidimensional case.

%%%%%%%%%%%%%%%%%%%%%%%%%%%%%
\section{Existence and uniqueness of solutions}\label{sec:WP}
In this section we construct solutions of the Fokker-Planck equation \eqref{eq1},
supplemented with the boundary condition
\( \label{BC1}
   D u_\balpha \grad_\bx \left(\frac{u}{u_\balpha}\right) \cdot \nu = 0 \quad\mbox{a.e. on } \partial\Omega_\bx,
\)
where $\nu=\nu(x)$ denotes the unit normal vector to the boundary of $\Omega_\bx$.
Moreover, we prescribe the initial condition
\(  \label{IC1}
   u(t=0) = u_0 \quad\mbox{on } \Omega_\bx.
\)
Our strategy is to convert the problem to the Hamiltonian form $(-\laplace + V)$
for a suitable potential $V$ and construct the corresponding semigroup.
In order to obtain some intuition, we first carry out the transform formally.

%%%%%%%%%%%%%%%%%%%%%%%%%
\subsection{Formal calculations}
Setting
\[
   \AA(\bx) := \sum_{i=1}^L \ln\xi(x_i) - 2N\aA(\bx),
\]
\eqref{eq1} is written in the form
\[
   \part{u}{t} = \grad_\bx\cdot\bigl( D(\grad_\bx u + u\grad_\bx\AA )\bigr)
\]
with the boundary condition $D(\grad_\bx u + u\grad_\bx\AA )\cdot\nu = 0$.
For $i=1,\dots,N$ we introduce the coordinate transform
\(   \label{xtoy}
   y_i := y(x_i) := 2\sqrt{N} \int_0^{x_i} \frac{\d s}{\sqrt{\xi(s)}} =  4\sqrt{N} \arcsin\sqrt{x_i},
\)
and denote $\by := (y_1,\dots,y_L)$. Note that $\bx\mapsto\by$ maps $\Omega_\bx = (0,1)^L$
onto $\Omega_\by := (0, Y_N)^L$ with $Y_N:=2 \pi\sqrt{N}$.
Introducing the new variable
\(   \label{baru}
   \bar u(\by) := J(\bx(\by)) u(\bx(\by)), \qquad J(\bx):= \left(2\sqrt{N}\right)^{-L} \prod_{j=1}^L \xi^{1/2}(x_i)
\)
transforms \eqref{eq1} to the form
\(   \label{eq1transf}
   \part{\bar u}{t} = \grad_\by\cdot\left( \grad_\by\bar u + \bar u \grad_\by\left( \bar\AA - \ln \bar J \right) \right),
\)
with $\bar\AA(\by) := \AA(\bx(\by))$ and $\bar J(\by) := J(\bx(\by))$.
By $\bx(\by)$ we denote the componentwise inverse transform $x_i=x(y_i)$.
The no-flux boundary condition \eqref{BC1} transforms as
\[
   \frac{1}{2\sqrt{N}}  \sum_{i=1}^L \bar J^{-1} \sqrt{\xi_i} \left( \partial_{y_i} \bar u + \bar u \partial_{y_i}
     \left( \bar\AA - \ln \bar J \right) \right) \nu_i = 0 \quad\mbox{a.e. on } \partial\Omega_\by,
\]
where we use the shorthand notation $\xi_i = \xi(x_i(y_i))$.
Note that the product $\bar J^{-1} \sqrt{\xi_i}$ is constant in $y_i$ and positive on the
set $\bigl\{\by\in\partial\Omega_\by; y_i\in\{0,Y_N \bigr\}, 0 < y_j < Y_N \mbox{ for } j\neq i\}$.
Consequently, the transformed boundary condition is equivalent to the nondegenerate expression
\(  \label{BCtransf}
    \left( \grad_\by \bar u + \bar u \grad_\by \left( \bar\AA - \ln \bar J \right) \right) \cdot \nu = 0 \quad\mbox{ a.e. on } \partial\Omega_\by,
\)
which can be also written as $\nu\cdot\grad_\by \ln (\bar u/\bar u_\balpha) = 0$ a.e. on $\partial\Omega_\by$.
The steady state for \eqref{eq1transf}--\eqref{BCtransf} is
\[
   \bar u_\balpha := \overline{\mathbb{Z}}_\balpha^{-1} \exp(-(\bar\AA - \ln \bar J)), \qquad
   \overline{\mathbb{Z}}_\balpha := \int_{\Omega_\by} \exp(-(\bar\AA - \ln \bar J)) \d\by.
\]
Finally, setting
\[
   z(\by):=\bar u(\by) / \sqrt{\bar u_\balpha(\by)},
\]
the Fokker-Planck equation \eqref{eq1} transforms to the Hamiltonian form
\(   \label{zeq}
   \part{z}{t} = \laplace_\by z - V(\by)z,
\)
with
\[
   V(\by) = \frac{\laplace_\by \sqrt{\bar u_\balpha}}{\sqrt{\bar u_\balpha}}
      = \frac12 \frac{\laplace_\by \bar u_\balpha}{\bar u_\balpha} - \frac14 \frac{|\grad_\by \bar u_\balpha|^2}{\bar u_\balpha^2},
\]
which can be further expressed as
\(  \label{V}
   V(\by) = -\frac12 \laplace_\by (\bar\AA - \ln \bar J) + \frac14 |\grad_\by (\bar\AA - \ln \bar J))|^2.
\)
The boundary condition \eqref{BCtransf} transforms to
\(  \label{zBC}
    \sqrt{\bar u_\balpha}  \left( \grad_\by z + \frac12 z \grad_\by \left( \bar\AA - \ln \bar J \right) \right)\cdot \nu = 0
    \qquad\mbox{a.e. on } \partial\Omega_\by.
\)
Let us remark that with \eqref{alphaA}, $\bar u_\balpha$ behaves like $\left(\prod_{i=1}^L \xi_i\right)^{4N\mu-1/2}$ close to the boundary,
so that for $4N\mu-1/2 > 0$ the boundary condition \eqref{zBC} is degenerate.

Inserting the expression \eqref{alphaA} for $\aA$ into \eqref{V}
gives the explicit expression for the potential
\(  \label{Vexp}
   V = \frac{1}{16N} \left(4N\mu - \frac12\right)\left(4N\mu - \frac32\right) \frac{(\xi'_i)^2}{\xi_i} + \mbox{(bounded terms)},
\)
where (bounded terms) are expressions involving
\[
   \part{\xi_i}{y_i} = \frac{\sqrt{\xi_i}}{2\sqrt{N}} \xi_i',\qquad
   \partk{\xi_i}{y_i}{2} =-\frac{1}{2N} \xi_i + \frac{1}{4N} (\xi_i')^2,\\
   \part{\xi'_i}{y_i} = - \frac{1}{\sqrt{N}} \sqrt{\xi_i},\qquad
   \partk{\xi'_i}{y_i}{2} = \frac{1}{4N} \xi'_i,
\]
that are uniformly bounded on $\overline\Omega_\by$.
The unbounded term in $V$ is
\[
    \frac{(\xi'_i)^2}{\xi_i} = \frac{(1-2x_i)^2}{x_i(1-x_i)},
\]
so for the potential to be bounded below, we need $4N\mu \geq 3/2$.

%%%%%%%%%%%%%%%%%%%%%%%%
\subsection{Construction of solutions for the case $4N\mu \geq 1/2$}
In this Section we shall construct weak solutions of the Fokker-Planck equation \eqref{eq1} with $4N\mu \geq 1/2$,
subject to the no-flux boundary condition \eqref{BC1} and the initial datum \eqref{IC1}.
However, since the equivalent form \eqref{zeq} is more suitable to study the asymptotic
behavior of the solution for large times, we shall work with this formulation.
Due to the issues caused by the degeneracy of the boundary condition,
we shall start from a weak formulation of \eqref{eq1} and carry out the coordinate transform
as in previous Section in order to arrive at a weak formulation of \eqref{zeq}.

To obtain a symmetric form, we multiply \eqref{eq1} by $\phi/u_\balpha$,
with a test function $\phi\in C^\infty(\overline\Omega_\bx)$,
and integrate by parts, taking into account the no-flux boundary condition \eqref{BC1}.
We arrive at
\[ \label{weak1}
   \tot{}{t} \int_{\Omega_\bx} \frac{u}{u_\balpha} \frac{\phi}{u_\balpha} u_\balpha \d \bx = - \int_{\Omega_\bx} D \grad_\bx \left(\frac{u}{u_\balpha}\right)
      \cdot  \grad_\bx \left(\frac{\phi}{u_\balpha}\right) u_\balpha \d \bx.
\]
\def\dbuy{\d\bar u_\balpha(\by)}
Carrying out the coordinate transform $\bx\mapsto\by$ \eqref{xtoy}, with the Jacobian $J$ given by
\eqref{baru}, yields
\[
   \tot{}{t} \int_{\Omega_\by} \frac{\bar u}{\bar u_\balpha} \frac{\bar \phi}{\bar u_\balpha} \dbuy =
     - \int_{\Omega_\by}  \grad_\by \left(\frac{\bar u}{\bar u_\balpha}\right)
      \cdot  \grad_\by \left(\frac{\bar\phi}{\bar u_\balpha}\right) \dbuy,
\]
with $\bar u$ given by \eqref{baru}, $\bar\phi(\by) := J(\bx(\by)) \phi(\bx(\by))$
and $\dbuy := \bar u_\balpha \d \by$.
Finally, defining $z:=\bar u/\sqrt{\bar u_\balpha}$ and $\psi:=\bar\phi/\sqrt{\bar u_\balpha}$,
we arrive at
\( \label{weakz}
   \tot{}{t} \int_{\Omega_\by} \frac{z}{\sqrt{\bar u_\balpha}} \frac{\psi}{\sqrt{\bar u_\balpha}} \dbuy =
     - \int_{\Omega_\by}  \grad_\by \left(\frac{z}{\sqrt{\bar u_\balpha}}\right)
      \cdot  \grad_\by \left(\frac{\psi}{\sqrt{\bar u_\balpha}}\right) \dbuy.
\)
We thus define the space
\( \label{H}
   \H_\by := \left\{ z\in L^2(\Omega_\by); \int_{\Omega_\by} \left| \grad_\by \frac{z}{\sqrt{\bar u_\balpha}} \right|^2 \dbuy < +\infty \right\}
\)
with the scalar product
\[
   (z, \psi)_{\H_\by} :=  \int_{\Omega_\by} \frac{z}{\sqrt{\bar u_\balpha}} \frac{\psi}{\sqrt{\bar u_\balpha}} \dbuy  +
      \int_{\Omega_\by}  \grad_\by \left(\frac{z}{\sqrt{\bar u_\balpha}}\right)
      \cdot  \grad_\by \left(\frac{\psi}{\sqrt{\bar u_\balpha}}\right) \dbuy
\]
and the induced norm $\Norm{z}_{\H_\by}^2 := (z,z)_{\H_\by}$.
Central for our analysis is the following result.

\begin{lemma}\label{lem:ineq1}
Let $4N\mu \geq 1/2$.
Then for every $z\in\H_\by$ the inequality holds
\[
   \int_{\Omega_\by} \left| \frac{z}{\sqrt{\bar u_\balpha}} \right|^2 \dbuy \geq \int_{\Omega_\by} |\grad z|^2 + Vz^2 \d\by,
\]
with $V$ defined in \eqref{V}.
\end{lemma}

\begproof
We have
\[
   \int_{\Omega_\by} \left| \frac{z}{\sqrt{\bar u_\balpha}} \right|^2 \dbuy =
      \int_{\Omega_\by} \left( |\grad z|^2 + \frac14 \frac{|\grad_\by \bar u_\balpha|^2}{\bar u_\balpha^2}z^2
      - \frac{\grad_\by z \cdot \grad_\by \bar u_\balpha}{\bar u_\balpha} z \right) \d\by.
\]
We integrate by parts in the last term of the right-hand side,
\[
   - \int_{\Omega_\by} \frac{\grad_\by z \cdot \grad_\by \bar u_\balpha}{\bar u_\balpha} z \d\by
   &=& -\frac12 \int_{\Omega_\by} \frac{\grad_\by \bar u_\balpha}{\bar u_\balpha} \cdot \grad_\by z^2 \d\by \\
   &=& \frac12 \int_{\Omega_\by} z^2 \grad_\by\cdot \left( \frac{\grad_\by \bar u_\balpha}{\bar u_\balpha} \right) \d\by
        - \frac12 \int_{\partial\Omega_\by} z^2 \frac{\grad_\by \bar u_\balpha}{\bar u_\balpha}\cdot\nu \d S_\by.
\]
With \eqref{alphaA} we have
\[
   \frac{\grad_{y_i} \bar u_\balpha}{\bar u_\balpha} = \grad_{y_i} (\ln\bar J - \bar\AA) = 
     \frac{\sqrt{\xi_i}}{2\sqrt{N}} \left[ \left(4N\mu-\frac12\right) \frac{\xi_i'}{\xi_i} - 4N\gamma_i + 4N\eta\xi_i' \right].
\]
Since $\xi_i$ vanishes for $y_i\in\{0, Y_N\}$ and $\xi_i'$ is bounded on $[0,Y_N]$, we have
\[
   - \frac12 \int_{\partial\Omega_\by} z^2 \frac{\grad_\by \bar u_\balpha}{\bar u_\balpha}\cdot\nu \d S_\by
     =  - \frac{1}{4\sqrt{N}} \left(4N\mu-\frac12\right) \sum_{i=1}^L \int_{\partial\Omega_\by}  z^2 \frac{\xi_i'}{\sqrt{\xi_i}} \nu_i \d S_\by.
\]
We write the boundary of the hypercube $\Omega_\by$ as an union of the pairs of faces,
\[
   \partial\Omega_\by = \bigcup_{i=1}^L F_i,\qquad F_i:=\{\by\in \partial\Omega_\by,\, y_j \in\{0,Y_N\}\},
\]
then we have
\[
    \sum_{i=1}^L \int_{\partial\Omega_\by}  z^2 \frac{\xi_i'}{\sqrt{\xi_i}} \nu_i \d S_\by =
        \sum_{i=1}^L \int_{F_i}  \left[ z^2 \frac{\xi_i'}{\sqrt{\xi_i}} \right]_{y_i=0}^{Y_N} \d S_{F_i},
\]
where $\d S_{F_i}$ denotes the $(L-1)$-dimensional Lebesgue measure on $F_i$.
Since $x_i'(x(Y_N)) = x_i'(1) = -1$ and $x_i'(x(0)) = x_i'(0) = 1$, we have
\[
   \left[ z^2 \frac{\xi_i'}{\sqrt{\xi_i}} \right]_{y_i=0}^{Y_N} \leq 0.
\]
Therefore, if $4N\mu-\frac12 \geq 0$,
\[
   - \frac12 \int_{\partial\Omega_\by} z^2 \frac{\grad_\by \bar u_\balpha}{\bar u_\balpha}\cdot\nu \d S_\by \geq 0.
\]
Consequently,
\[
    \int_{\Omega_\by} \left| \frac{z}{\sqrt{\bar u_\balpha}} \right|^2 \dbuy &\geq& 
    \int_{\Omega_\by} \left( |\grad z|^2 + \frac14 \frac{|\grad_\by \bar u_\balpha|^2}{\bar u_\balpha^2}z^2
      + \frac12  z^2 \grad_\by\cdot \left( \frac{\grad_\by \bar u_\balpha}{\bar u_\balpha} \right) \right) \d\by \\
      &=&
      \int_{\Omega_\by}  |\grad z|^2 + Vz^2 \d\by.
\]
Finally, the above formal calculation are made rigorous by replacing $\bar u_\balpha$
by $\bar u_\balpha^\eps := \bar u_\balpha + \eps$ for $\eps>0$ and subsequently passing
to the limit $\eps\to 0$.
\endproof

\begin{lemma}
Let $4N\mu \geq 1/2$.
Then the space $\H_\by$ defined in \eqref{H} with the scalar product $(\cdot,\cdot)_{\H_\by}$ is a Hilbert space,
and is densely embedded into $L^2(\Omega_\by)$.
\end{lemma}

\begproof
Completeness follows from the fact that if $z_k$ is a Cauchy sequence
in $\H_\by$, then due to Lemma \ref{lem:ineq1} it is also a Cauchy sequence in $L^2$.
The density of the embedding into $L^2(\Omega_\by)$ is due to the fact that the set
of smooth functions with compact support is dense in $\H_\by$.
\endproof

\begin{defin}\label{defin:weak}
We call $z\in L^2((0,T); \H_\by) \cap C([0,T]; L^2(\Omega_\by))$ a weak solution of \eqref{zeq} on $[0,T)$
subject to the boundary condition \eqref{zBC} if \eqref{weakz} holds for every $\psi\in \H_\by$
and almost all $t\in(0,T)$, and the initial condition is satisfied by continuity in $C([0,T]; L^2(\Omega_\by))$.
\end{defin}

We remark that a formal integration by parts in the right-hand side of \eqref{weakz} gives
\[
   - \int_{\Omega_\by}  \grad_\by \left(\frac{z}{\sqrt{\bar u_\balpha}}\right)
      \cdot  \grad_\by \left(\frac{\psi}{\sqrt{\bar u_\balpha}}\right) \dbuy \\
      = - \int_{\partial \Omega_\by} \left[ \sqrt{\bar u_\balpha} \left( \grad_\by z + \frac12 z \grad_\by (\bar\AA-\ln \bar J)\right)\cdot\nu \right]
      \frac{\psi}{\sqrt{\bar u_\balpha}} \d S_\by  \\
      + \int_{\Omega_\by} [ \laplace_\by z - V z ] \frac{\psi}{\sqrt{\bar u_\balpha}}  \d\by.
\]
This justifies the interpretation of \eqref{weakz} as the weak formulation of  \eqref{zeq} subject
to the boundary condition \eqref{zBC}.

We now define the operator $\L: D(\L) \subset \H_\by \to L^2(\Omega_\by)$ by its action
\(  \label{L}
   \langle \L z, \psi \rangle := - \int_{\Omega_\by}  \grad_\by \left(\frac{z}{\sqrt{\bar u_\balpha}}\right)
      \cdot  \grad_\by \left(\frac{\psi}{\sqrt{\bar u_\balpha}}\right) \dbuy
\)
for all $z$, $\psi\in \H_\by$.
We shall prove that the closure $\overline \L$ of $\L$ generates a contraction semigroup on $L^2(\Omega_\by)$.
For this sake, we study the resolvent problem
\(  \label{resolventP}
   (-\L + \lambda) z = f
\)
for (some) $\lambda>0$ and $f \in L^2(\Omega_\by)$.

\begin{lemma}\label{lem:res}
Let $4N\mu \geq 1/2$.
Then for every $f\in L^2(\Omega_\by)$ the resolvent problem \eqref{resolventP}
has a unique solution $z\in \H_\by$. 
\end{lemma}

\begproof
For a fixed $\lambda>0$ we define the bilinear form $a: \H_\by\times \H_\by \to \R$,
\[
   a_\lambda(z,\psi) := \langle -\L z,\psi \rangle + \lambda (z,\psi),
\]
where $(z,\psi)$ denotes the standard scalar product on $L^2(\Omega_\by)$.
The resolvent problem \eqref{resolventP} with the no-flux boundary conditions is equivalent to
\[
   a_\lambda (z,\psi) = (f, \psi) \qquad\mbox{for all } \psi\in \H_\by.
\]
A straightforward application of the H\"older inequality gives the continuity of $a_\lambda$,
\[
  a_\lambda(z,\psi) \leq C \Norm{z}_{\H_\by}\Norm{\psi}_{\H_\by}
\]
for a suitable constant $C>0$; coercivity is straightforward.
Finally, the mapping $\psi \mapsto (f,\psi)$ with $f\in L^2(\Omega_\by)$ is an element of the 
dual space $(\H_\by)'$. Consequently, an application of the Lax-Milgram theorem
yields the existence and uniqueness of the solution $z\in \H_\by$.
\endproof

\begin{theorem}\label{thm:semigroup}
Let $4N\mu \geq 1/2$.
Then the closure $\overline \L$ of $\L$ generates a contraction semigroup on $L^2(\Omega_\by)$. 
\end{theorem}

\begproof
Since $\H_\by$ is densely embedded into $L^2(\Omega_\by)$, the operator $\L$ is densely defined, and dissipative.
Moreover, due to Lemma \ref{lem:res}, the range of $-\L + \lambda$ is $L^2(\Omega_\by)$ for all $\lambda > 0$.
The claim then follows by an application of the Lumer-Phillips theorem \cite{Pazy}.
\endproof

The contraction semigroup constructed in Theorem \ref{thm:semigroup}
provides the announced existence and uniqueness of weak solutions
$z\in L^2((0,T); \H_\by) \cap C([0,T]; L^2(\Omega_\by))$ of \eqref{zeq}
subject to the no-flux boundary condition \eqref{zBC}
in the sense of Definition \ref{defin:weak}.
The solutions are formally written as $z(t) = e^{\L t}z_0$, where $z_0\in L^2(\Omega_\by)$
is the initial datum; see, e.g., \cite{Pazy}.
By the inverse coordinate transform to \eqref{xtoy}
we obtain weak solutions of the original Fokker-Planck equation \eqref{eq0}
subject to the no-flux boundary condition \eqref{BC1}.

%%%%%%%%%%%%%%%%%%%%%%%%%%%%%%%%%%%%
\section{Spectral gap - exponential convergence to equilibrium}\label{sec:gap}
In this Section we shall perform a spectral analysis of the operator $(-\L)$
and prove that boundedness below of the potential $V$ \eqref{V} implies
exponential convergence to equilibrium for \eqref{zeq}.
From the explicit expression \eqref{Vexp} for $V$ we see that
$V$ is bounded below if $4N\mu \geq 3/2$.

\begin{lemma}\label{lem:compact_resolvent}
Let $4N\mu \geq 3/2$.
Then the operator $(-\L)$ defined in \eqref{L} has compact resolvent.
\end{lemma}

\begproof
We need to show that for some $\lambda>0$ the operator $(-\L + \lambda)^{-1}$ is compact
as a mapping from $L^2(\Omega_\by)$ into itself.
Let $f\in L^2(\Omega_\by)$ and $z=(-\L + \lambda)^{-1} f$, constructed in Lemma \ref{lem:res}.
From Lemma \ref{lem:ineq1} we have
\[
   ((-\L + \lambda)z, z) \geq \int_{\Omega_\by} |\grad_\by z|^2 + (V + \lambda) z^2 \d\by \geq C \Norm{z}^2_{H^1(\Omega_\by)}
\]
for some constant $C>0$ and $\lambda$ chosen such that $\min_{\by\in\Omega_\by} (V(\by)+\lambda) > 0$.
On the other hand, the Cauchy-Schwartz inequality gives
\[
  ((-\L + \lambda)z, z) = (f,z) \leq \frac{1}{2\eps} \Norm{f}_{L^2(\Omega_\by)}^2 + \frac{\eps}{2} \Norm{z}_{L^2(\Omega_\by)}^2,
\]
so for sufficiently small $\eps>0$ we conclude
\[
    \Norm{(-\L + \lambda)^{-1} f}_{H^1(\Omega_\by)} = \Norm{z}_{H^1(\Omega_\by)} \leq C \Norm{f}_{L^2(\Omega_\by)}
\]
and the claim follows by the compact embedding of the Sobolev space $H^1$ into $L^2$.
\endproof

Together with the obvious self-adjointness of $(-\L)$, Lemma \ref{lem:compact_resolvent} implies that  $(-\L)$
has a discrete spectrum without finite accumulation points. Moreover, all its eigenvalues are nonnegative.
This implies the existence of a positive spectral gap and, consequently, exponential convergence to equilibrium
as $t\to\infty$, see, e.g., \cite{AMTU}.

%%%%%%%%%%%%%%%%%%
\section{The Dynamical Maximum Entropy Approximation}\label{sec:DNE}
In typical applications in quantitative genetics the solution of the
Fokker-Planck equation \eqref{eq0} is not the main object of interest.
One is rather interested in the evolution of its certain moments
that correspond to the macroscopic dynamics of observable quantitative traits.
This naturally leads to the question whether one can derive a finite-dimensional
system of differential equations that approximates the evolution
of the moments of interest, avoiding the need of solving \eqref{eq0}.
This question has been studied previously by analogy with statistical mechanics:
the allele frequency distribution is approximated by the stationary
form, which maximizes the logarithmic relative entropy.
Called Maximum Entropy Method,
it has been applied to broad spectrum of problems ranging from the statistics of neural spiking
\cite{Schneidman, Tkacik}, bird flocking
\cite{Bialek}, protein structure \cite{Weigt}, immunology \cite{Mora} and more.
For transient problems described by known dynamical equations (e.g., Fokker-Planck equation),
the \emph{Dynamical Maximum Entropy (DynMaxEnt)}
method assumes quasi-stationarity at each time point.
It has been applied, e.g., to modeling of cosmic ray transport \cite{Hick},
general Fokker-Planck equation \cite{Plastino},
analysis of genetic algorithms \cite{Shapiro97},
and population genetics \cite{Shapiro01, Barton09, BTB}.
In \cite{BTB} it is observed that the "classical" DynMaxEnt method cannot be applied
in the regime of small mutations, and the theory is extended for this regime
to account for changes in mutation strength.
Surprisingly, systematic numerical simulations document
superb approximation properties of the method even far from the quasi-stationary regime.
However, derivation of analytic error estimates remains an open problem.
%The accuracy of the , however, remain unexplained.

In this section we discuss several aspects of the DynMaxEnt method.
First, in Section \ref{subsec:entropy} we show that constrained maximization of a logarithmic entropy
functional leads to a moment-matching condition.
%which we solve in the simple case \ref{ssec:solvability}.
Then, in Section \ref{subsec:DNE} we provide a simple and straightforward derivation
of the DynMaxEnt method by adopting a quasi-stationary approximation.
To our best knowledge, this derivation has not been known before.
Finally, in Section \ref{subsec:simple} we consider the scalar case and derive a modified version of
the DynMaxEnt method, which is valid for the whole range of admissible parameters.

%%%%%%%%%%%%%%%%%%
\subsection{Constrained entropy maximization}\label{subsec:entropy}
We shall call the vector $\balpha\in\R^d$ \emph{admissible} if the corresponding
normalization factor $\mathbb{Z}_\balpha$ \eqref{Z_alpha}
is finite.
For any integrable function $u\in L^1(\Omega_\bx)$ with $\int_{\Omega_\bx} u(\bx) \d\bx = 1$
and any admissible $\balpha\in\R^d$ we define the logarithmic relative entropy
\( \label{relEnt}
   H(u | u_\balpha) := \int_{\Omega_\bx} u\ln\frac{u}{u_\balpha} \d\bx,
\)
where $u_\balpha$ is the normalized stationary solution of the Fokker-Planck equation \eqref{eq0},
given by formula \eqref{stat0}. Note that this is a different approach compared with \cite{BTB},
where the logarithmic entropy is taken relative to the neutral distribution of allele frequencies in
the absence of mutation or selection, $\prod_{i=1}^L \xi_i^{-1}$,
and the variational problem is complemented with normalization and moment constraints.

For a fixed $u\in L^1(\Omega_\bx)$ with finite $\bA$-moments, let us consider the maximization of the relative entropy \eqref{relEnt}
in terms of admissible $\balpha\in\R^d$, i.e., the task of maximizing the function $\balpha \mapsto H(u | u_\balpha)$.
If a critical point exists, then for $i=1,\dots,d$,
\[
   \tot{}{\alpha_i} H(u | u_\balpha) &=& - \int_{\Omega_\bx} \frac{u}{u_\balpha} \tot{}{\alpha_i} u_\balpha \d\bx \\
    &=& \mmnt{A_i}_{u_\balpha} - \mmnt{A_i}_u = 0.
\]
Consequently, if a maximizer $\balpha^*$ exists, then the $\bA$-moments corresponding to $u_{\balpha^*}$
must be matching the same moments of $u$.
This naturally leads to the question of solvability of the nonlinear
system of equations
\[
   \mmnt{\bA}_{u_\balpha}  =  \mmnt{\bA}_u
\]
in terms of the admissible parameter vector $\balpha\in\R^d$, for a given, normalized $u\in L^1(\Omega_\bx)$
with finite $\bA$-moments.
To address this question seems to be a very difficult task that we leave open.
We merely remark that the Hessian matrix of $\balpha \mapsto H(u | u_\balpha)$,
\[
    \frac{\d^2}{\d\alpha_i \d\alpha_j} H(u | u_\balpha) = \mmnt{A_i A_j}_{u_\balpha} - \mmnt{A_i}_{u_\balpha} \mmnt{A_j}_{u_\balpha},
\]
is equal to the covariance matrix of the random variables $\bA$ with the probability density $u_\balpha$.
Thus, the Hessian matrix is positive semidefinite.
In the scalar case, solvability of the moment equation $ \mmnt{A}_{u_\alpha}  =  \mmnt{A}_u$
can be studied for particular choices of $A$. We shall give an example below
in Section \ref{ssec:solvability}.

%%%%%%%%%%%%%%%%%%
\subsection{Derivation of the DynMaxEnt method}\label{subsec:DNE}
Let us consider $u=u(t)$ a solution of the Fokker-Planck equation \eqref{eq0} with
admissible parameter vector $\balpha$,
subject to the initial datum $u(t=0) = u_{\balpha^0}$ for some admissible $\balpha^0$.
The DynMaxEnt method is derived in two steps:
First, we multiply the equation in its form \eqref{eq1} by the vector $\bA$
and integrate,
\[
   \tot{}{t} \mmnt{\bA}_{u(t)} &=& \int_{\Omega_\bx} \bA \grad_\bx \left( D u_\balpha \grad_\bx \left( \frac{u}{u_\balpha} \right) \right) \d\bx \\
     &=& - \int_{\Omega_\bx} \grad_\bx \bA  D u_\balpha \grad_\bx \left( \frac{u}{u_\balpha} \right) \d\bx,
\]
where we assumed that the boundary term in the integration by parts vanishes
(note that, in general, this does not necessarily follow from \eqref{BC1}).
In the second step, we substitute $u(t)$ in the above expression by $u_{\balpha^*(t)}$
with some time-dependent parameter vector $\balpha^*=\balpha^*(t)$,
which leads to
\[
   \tot{}{t} \mmnt{\bA}_{u_{\balpha^*(t)}} = - \int_{\Omega_\bx} \grad_\bx \bA  D u_\balpha \grad_\bx \left( \frac{u_{\balpha^*(t)}}{u_\balpha} \right) \d\bx + \mathbf{R},
\]
where $\mathbf{R}$ is a vector-valued residuum term. We now introduce an approximation by neglecting the residuum $\mathbf{R}$.
Expanding the derivatives on both sides of the above equation leads then to
\(  \label{DynMaxEnt}
 \left( \mmnt{\bA\otimes\bA}_{u_{\balpha^*(t)}}
     - \mmnt{\bA}_{u_{\balpha^*(t)}} \otimes \mmnt{\bA}_{u_{\balpha^*(t)}} \right) \tot{\balpha^*(t)}{t}
    =  \\
    = \frac12 \mmnt{\xi \grad_\bx \bA : \grad_\bx \bA}_{u_{\balpha^*(t)}} (\balpha - \balpha^*(t)),
    \nonumber
\)
where $\grad_\bx \bA : \grad_\bx \bA$ is the symmetric $d\times d$ matrix with the $(i,k)$-component
$\sum_{j=1}^d \partial_{x_j} A_i \partial_{x_j} A_k$.
The nonlinear ODE system for $\balpha^*=\balpha^*(t)$ is called the \emph{DynMaxEnt method} for approximation
of the moments of \eqref{eq0}.
However, two comments have to be made:
First, the matrix on the left-hand side,% multiplying the time-derivative of $\balpha^*(t)$
$$\left( \mmnt{\bA\otimes\bA}_{u_{\balpha^*(t)}}
     - \mmnt{\bA}_{u_{\balpha^*(t)}} \otimes \mmnt{\bA}_{u_{\balpha^*(t)}} \right),$$
is positive semidefinite, since it is the covariance matrix of the observables $\bA$
of the probability distribution $u_{\balpha^*(t)}$. However, in order
\eqref{DynMaxEnt} to be globally solvable, the covariance matrix must be uniformly (positive) definite,
which in general may not be the case.
Furthermore, the matrix $\mmnt{\xi \grad_\bx \bA : \grad_\bx \bA}_{u_{\balpha^*(t)}}$
may have infinite entries even for some admissible ${\balpha^*(t)}$, and if this is the case,
then again the ODE system is not solvable.
Since these two issues are very hard to resolve in general,
we shall below resort to a simple case where $\balpha$ is a scalar.

\subsection{Scalar case}\label{subsec:simple}
To gain some more insight into the ODE \eqref{DynMaxEnt},
we consider the single locus case $x\in(0,1)$
with $\bA$ being a scalar function $A=A(x)$
and $\alpha\in\R$.
The DynMaxEnt method \eqref{DynMaxEnt} simplifies to
the following ODE for $\alpha^* = \alpha^*(t)$,
\(  \label{DNEs}
    \left( \mmnt{A^2}_{u_{\alpha^*(t)}}
     - \mmnt{A}^2_{u_{\alpha^*(t)}} \right) \tot{\alpha^*(t)}{t}
    =  \frac12 \mmnt{\xi (\partial_x A)^2}_{u_{\alpha^*(t)}} (\alpha - \alpha^*(t)).
\)
An application of the Cauchy-Schwartz inequality implies that
\[
   \mmnt{A^2}_{u_{\alpha^*(t)}} - \mmnt{A}^2_{u_{\alpha^*(t)}} \geq 0,
\]
and, moreover, equality holds if and only if $A$ is a constant function.
Consequently, for every nonconstant $A$ the ODE \eqref{DNEs}
can be rewritten as
\(  \label{DNEscalar}
   \tot{\alpha^*(t)}{t} =  \frac12 \left( \mmnt{A^2}_{u_{\alpha^*(t)}}
     - \mmnt{A}^2_{u_{\alpha^*(t)}} \right)^{-1}
   \mmnt{\xi (\partial_x A)^2}_{u_{\alpha^*(t)}} (\alpha - \alpha^*(t)).
\)
However, the question of finiteness of the moment $\mmnt{\xi (\partial_x A)^2}_{u_{\alpha^*(t)}}$
can be only answered by making a particular choice for $A=A(x)$.

As a toy model, let us choose $A=A(x)$ to be the scalar function $\ln(\xi(x))$.
This corresponds to a population of individuals in a neutral environment ($\beta = h = 0$ in \eqref{alphaA})
with the nonzero mutation rate $\alpha = 2\mu$.
With the singularities at $x\in\{0, 1\}$, the function $A(x)=\ln(\xi(x))$ well represents the issues
that one encounters with the generic choice \eqref{alphaA}.
It is easily checked that the set of admissible values of $\alpha$
is the interval $(0,\infty)$.
Moreover, the moment $\mmnt{\xi (\partial_x A)^2}_{u_{\alpha^*}}$
is only finite for $\alpha^* > 1$. Consequently, the DynMaxEnt method \eqref{DNEscalar} is only applicable
if both the initial value $\alpha^*(0) = \alpha^0$ and $\alpha$ are strictly larger than $1$.
Then, since obviously the solution $\alpha^*(t)$ of \eqref{DNEscalar} is a monotone function of time,
it will stay strictly larger than $1$ for all $t\geq 0$ and asymptotically converge to $\alpha$.

The issue of non-finiteness of the term $\mmnt{\xi (\partial_x A)^2}_{u_{\alpha^*(t)}}$
%for $A(x) := \ln(\xi(x))$ and $\alpha^* \leq 1$
was addressed in \cite{BTB} by introducing a special treatment near the boundary
(see Appendix E, equations E.10-E.13 of \cite{BTB} for details of the derivation of the modified method).
Here we propose an alternative way that treats the problem at least in the
case $A(x) := \ln(\xi(x))$. It is based on the idea of multiplying the Fokker-Planck equation
by a suitable function $B=B(x)$, instead of $A=A(x)$, and integrating on $\Omega_\bx$.
In the second step, one again approximates $u(t)$ by $u_{\balpha^*(t)}$ and neglects the residuum.
This leads, in the scalar case, to the ODE
\(  \label{DNEsB}
    \left( \mmnt{AB}_{u_{\alpha^*(t)}}
     - \mmnt{A}_{u_{\alpha^*(t)}} \mmnt{B}_{u_{\alpha^*(t)}} \right) \tot{\alpha^*(t)}{t}
    =  \frac12 \mmnt{\xi (\partial_x B)^2}_{u_{\alpha^*(t)}} (\alpha - \alpha^*(t)).
\)
Choosing $B(x) := \xi(x)$ leads then to finite $\mmnt{\xi (\partial_x B)^2}_{u_{\alpha^*}}$
for all $\alpha^* > 0$, i.e., for all admissible values of $\alpha^*$. % for which the equilibrium state $u_{\alpha^*}$ is integrable.
Thus, our strategy is to obtain $\alpha^*(t)$ by solving
\eqref{DNEsB} for $t\geq 0$ and then calculate the moment $\mmnt{\ln(\xi)}_{u_{\alpha^*(t)}}$,
which is expected to be a good approximation of the true moment $\mmnt{\ln(\xi)}_{u(t)}$.
Clearly, one can use this strategy to obtain an approximation of any other moment of $u(t)$.

\begin{figure}
{\centering
\resizebox*{0.65\linewidth}{!}{\includegraphics{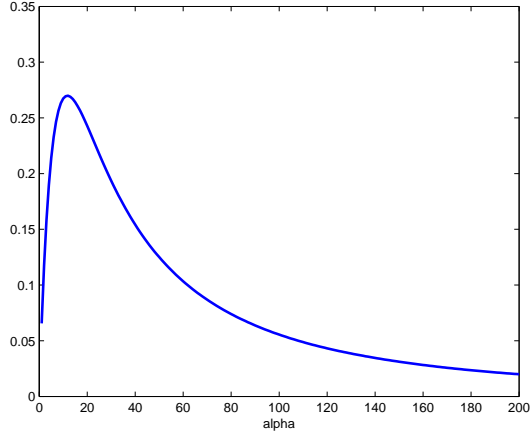}}
\par}
\caption{Numerical calculation of the expression $\mmnt{\xi\ln(\xi)}_{u_{\alpha^*}} - \mmnt{\ln(\xi)}_{u_{\alpha^*}} \mmnt{\xi}_{u_{\alpha^*}}$
for $\alpha^*\in [0.01, 200]$. Compare with the analytical result \eqref{limits}.
\label{fig:alpha_plot}}
\end{figure}

However, the method \eqref{DNEsB} suffers from a serious drawback, namely,
it is only solvable if the covariance
$$ \mmnt{AB}_{u_{\alpha^*(t)}} - \mmnt{A}_{u_{\alpha^*(t)}} \mmnt{B}_{u_{\alpha^*(t)}} $$
is nonvanishing for all $t\geq 0$, which is not clear.
Nonetheless, for the particular choice $A(x) = \ln(\xi(x))$ and $B(x) = \xi(x)$
this seems to be the case, as is documented by our numerical calculation in Fig. \ref{fig:alpha_plot}.
Analytically we are only able to calculate the limits
\( \label{limits}
    \mmnt{\xi\ln(\xi)}_{u_{\alpha^*}} - \mmnt{\ln(\xi)}_{u_{\alpha^*}} \mmnt{\xi}_{u_{\alpha^*}}
     \to 0 \qquad\mbox{as } \alpha^*\to 0, +\infty,
\)
which is based on the following Lemma.

\begin{lemma}\label{lem:limits}
For $\sigma>0$ and $x\in(0,1)$ denote
\(   \label{nu_sigma}
   \nu_\sigma(x) := \frac{\xi(x)^{\sigma-1}}{\int_0^1 \xi(s)^{\sigma-1} \d s}.
\)
with $\xi(x) = x(1-x)$.
Then, in the sense of distributions,
\[
   \nu_\sigma &\to& \delta(\cdot - 1/2) \qquad\mbox{as } \sigma\to \infty,\\ 
   \nu_\sigma &\to& \frac12 \delta(\cdot) + \frac12 \delta(\cdot - 1) \qquad\mbox{as } \sigma\to 0,\\ 
\]
where $\delta(\cdot - x)$ denotes the Dirac-delta distribution concentrated at $x$.
\end{lemma}

\begproof
Obviously, $\nu_\sigma$ is a probability measure on $(0,1)$.
Let $\phi\in C_c^\infty(0,1)$ be any test function on the interval $(0,1)$.
We shall show that
\[
   \lim_{\sigma\to\infty} \int_0^1 \phi(x) \d\nu_\sigma(x) = \phi(1/2).
\]
The mean-value theorem gives
\[
   \left| \phi(1/2) - \int_0^1 \phi(x) \d\nu_\sigma(x) \right| \leq \int_0^1 \left| \phi'(\eta(x)) \right| |x-1/2| \d\nu_\sigma(x)
     \leq C_\phi \int_0^1 |x-1/2| \d\nu_\sigma(x).
\]
Thus, our goal is to show that $\int_0^1 |x-1/2| \d\nu_\sigma(x)$ vanishes as $\sigma\to\infty$.
For the numerator, we have
\[
   \int_0^1 (4\xi(x))^{\sigma-1} |x-1/2| \d x = 2 \int_0^{1/2} (4\xi(x))^{\sigma-1} (1/2 - x) \d x,
\]
and using the identity $1/2 - x = \xi'(x)/2$, we calculate
\[
   \int_0^1 (4\xi(x))^{\sigma-1} |x-1/2| \d x = \frac{1}{\sigma}.
\]
The denominator is estimated from below using the elementary inequalities
\[
   4\xi(x) &\geq& 3x\qquad\mbox{for } x\in[0,1/4], \\ 
      &\geq& x+1/2 \qquad\mbox{for } x\in[1/4,1/2], \\ 
\]
which give
\[
   \int_0^1 (4\xi(x))^{\sigma-1} \d x \geq \frac{2}{\sigma} \left[ \left(\frac32\right)^\sigma - \frac23 \left(\frac34\right)^\sigma \right].
\]
Thus,
\[
   \int_0^1 |x-1/2| \d\nu_\sigma(x) \leq \frac12 \left[ \left(\frac32\right)^\sigma - \frac23 \left(\frac34\right)^\sigma \right]^{-1} \to 0
   \qquad\mbox{as } \sigma\to\infty,
\]
which proves the first claim.

To calculate the limit $\sigma\to 0$, due to the symmetry of $\xi(x)=x(1-x)$ with respect to $x=1/2$,
it is sufficient to prove that
\[
   \frac{\xi(x)^{\sigma-1}}{\int_0^{1/2} \xi(s)^{\sigma-1} \d s} \to \delta(\cdot) \qquad\mbox{as } \sigma\to 0.
\]
Again, picking a test function $\phi\in C_c^\infty[0,1/2)$
and using the mean-value theorem, we have to show that
\[
   \left| \phi(0) - \frac{\int_0^{1/2} \phi(x) \xi(x)^{\sigma-1}\d x}{\int_0^{1/2} \xi(s)^{\sigma-1} \d s} \right|
      \leq C_\phi \frac{\int_0^{1/2} \left| \phi'(\eta(x)) \right| x \xi(x)^{\sigma-1} \d x}{\int_0^{1/2} \xi(s)^{\sigma-1} \d s}
\]
tends to zero as $\sigma\to 0$. 
However, this follows directly from the fact that the numerator is uniformly bounded for, say, $0 \leq \sigma < 1$,
and that, obviously, the denominator tends to $+\infty$ as $\sigma\to 0$.
\endproof

The statement \eqref{limits} follows directly from the fact that
for $A(x) = \ln(\xi(x))$ we readily have $u_{\alpha^*} = \nu_{\alpha^*}$ with $\nu_{\alpha^*}$ given by \eqref{nu_sigma}.

Consequently, the "modified" DynMaxEnt method \eqref{DNEsB} can be safely used
with $A(x) = \ln(\xi(x))$ and $B(x) = \xi(x)$. It even seems to provide better approximation
results than the "original" method \eqref{DNEs}, as is documented by our numerical experiments
in Section \ref{sec:numerics}.

\subsubsection{Solvability of the moment equation}\label{ssec:solvability}
Finally, we study the solvability with respect to $\alpha>0$ of the moment equation
\(   \label{mmnt_eq}
   \mmnt{A}_{u_\alpha}  =  \mmnt{A}_u
\)
with $A(x) = \ln(\xi(x))$, assuming that the right-hand side is finite.
First of all, we note that the mapping $\alpha \mapsto \mmnt{A}_{u_\alpha}$
is strictly increasing for $\alpha > 0$. Indeed,
\[
   \tot{}{\alpha} \mmnt{A}_{u_\alpha} = \mmnt{A^2}_{u_\alpha} - \mmnt{A}^2_{u_\alpha} > 0,
\]
where the strict positivity follows as before by the Cauchy-Schwartz inequality.
Consequently, if a solution to the moment equation \eqref{mmnt_eq} exists, it is unique.
Next we claim that for any $u\geq 0$ with $\int_0^1 u(x) \d x = 1$ we have
$\mmnt{\ln(\xi)}_{u} \in [-\infty,\ln(1/4))$. Indeed, since $\xi(x) < 1/4$ on $(0,1)\setminus \{1/2\}$,
\[
   \mmnt{\ln(\xi)}_{u} < \ln\left( \frac14 \right) \int_0^1 u(x) \d x = \ln\left( \frac14 \right).
\]
Thus, it remains to prove that the range of $\alpha \mapsto \mmnt{A}_{u_\alpha}$
is the interval $(-\infty,\ln(1/4))$.
Since for $A(x) = \ln(\xi(x))$ we readily have $u_\alpha = \nu_\alpha$ with $\nu_\alpha$ given by \eqref{nu_sigma},
Lemma \ref{lem:limits} gives
\[
   \mmnt{\ln(\xi)}_{u_\alpha} = \int_0^1 \ln(\xi(x)) \d\nu_\alpha(x) &\to& \ln(1/4) \qquad\mbox{as } \alpha\to+\infty,\\
      &\to& -\infty \qquad\mbox{as } \alpha\to 0+.
\]
Indeed, the range of the mapping $\alpha \mapsto \mmnt{\ln(\xi)}_{u_\alpha}$ is the interval $(-\infty,\ln(1/4))$
and, therefore, the moment equation \eqref{mmnt_eq} is uniquely solvable for every normalized $u\in L^1(0,1)$
with finite $\mmnt{\ln(\xi)}_u$-moment.

\section{Numerical experiments}\label{sec:numerics}

\subsection{Scalar case}\label{subsec:num1}
We present results of several numerical experiments
that aim to demonstrate the performance of the
original \eqref{DNEs} and modified \eqref{DNEsB}
DynMaxEnt methods for the scalar (single locus) case $A(x) = \ln(\xi(x))$,
as discussed in Section \ref{subsec:simple}.
Let us recall that this case corresponds to a population of individuals
in a neutral environment ($\beta = h = 0$ in \eqref{alphaA})
with the nonzero mutation rate $\alpha = 2N\mu$.
For the modified method \eqref{DNEsB} we again choose $B(x) = \xi(x)$.

In all simulations we set $N=1$ and start from the initial condition $\alpha^*(t=0) = \alpha^0 := 2$
for the ODEs \eqref{DNEs}, \eqref{DNEsB}, and the initial datum
$u(t=0) = u_{\alpha^0}$ for the Fokker-Planck equation \eqref{eq1}.
The ODEs \eqref{DNEs}, \eqref{DNEsB} are solved with simple forward
Euler discretization on the time interval $[0,T]$ for different values of $T>0$.
We use $B(x) = \xi(x)$ for the modified DynMaxEnt method \eqref{DNEsB}.
The Fokker-Planck equation is discretized in space using the Chang-Cooper
scheme \cite{Chang-Cooper} and forward Euler in time.

\begin{figure}
{\centering \begin{tabular}[h]{cc}
\resizebox*{0.49\linewidth}{!}{\includegraphics{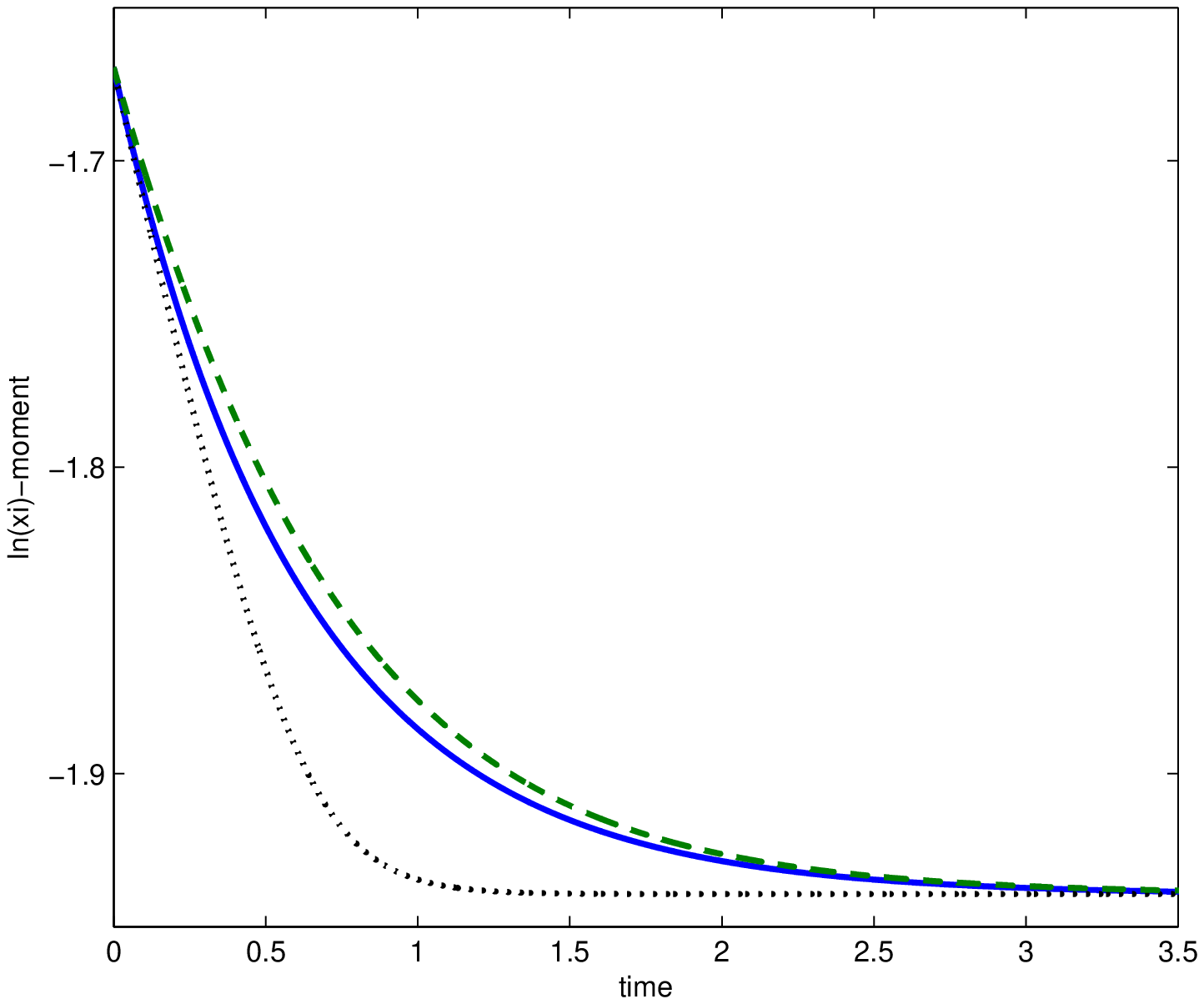}} &
\resizebox*{0.49\linewidth}{!}{\includegraphics{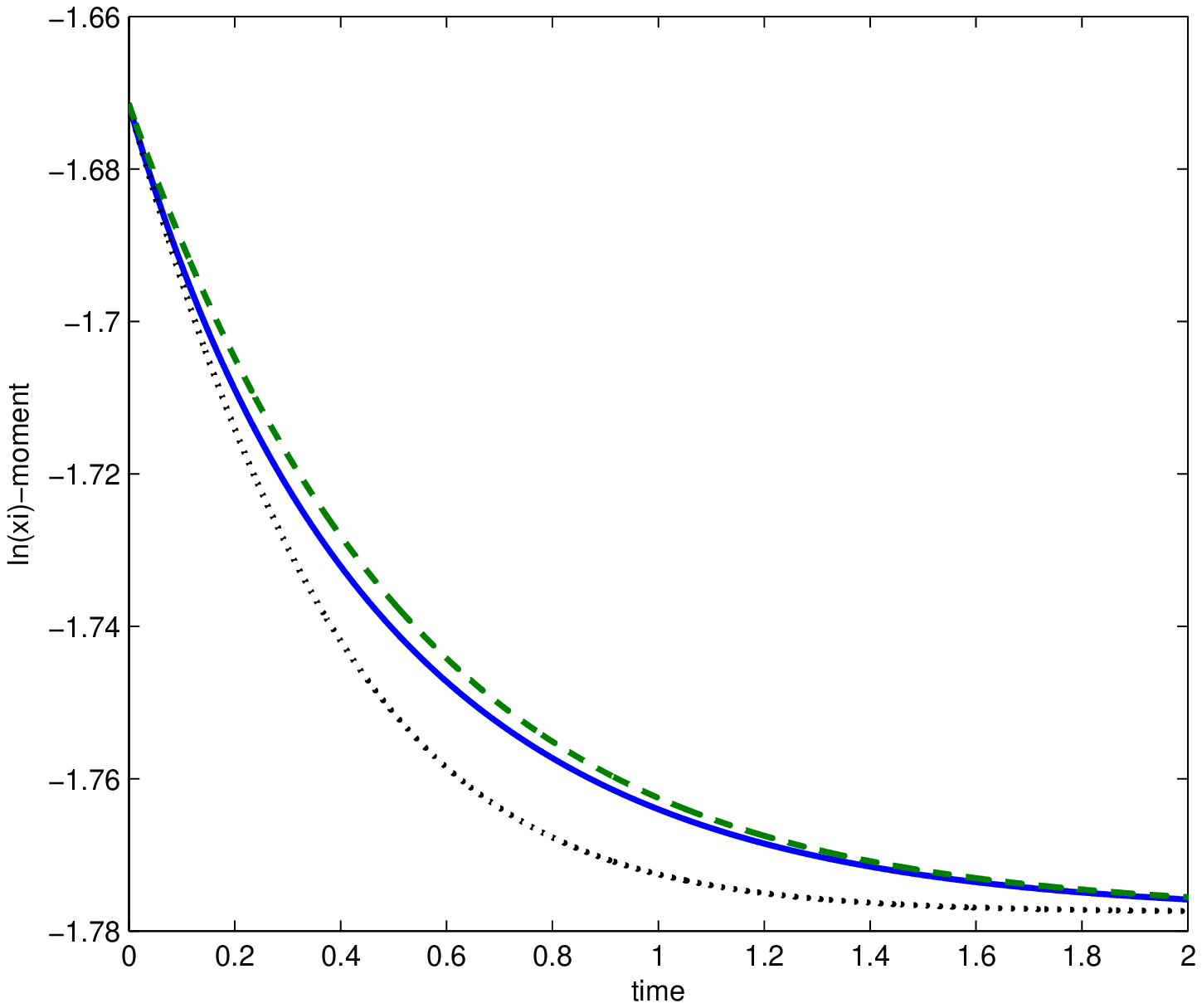}} \\
$\alpha = 1.1$ & $\alpha = 1.5$ \\
\resizebox*{0.49\linewidth}{!}{\includegraphics{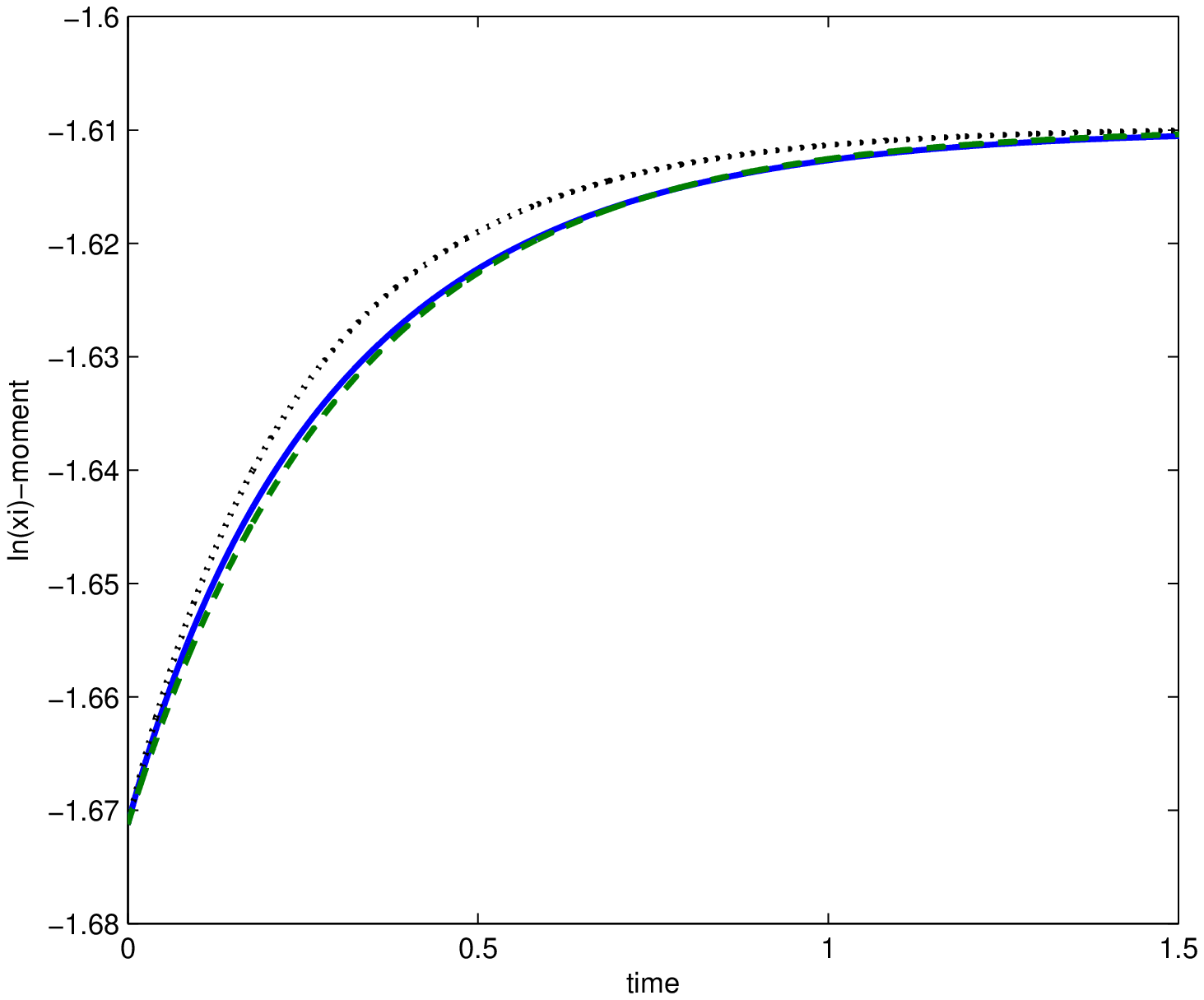}} &
\resizebox*{0.49\linewidth}{!}{\includegraphics{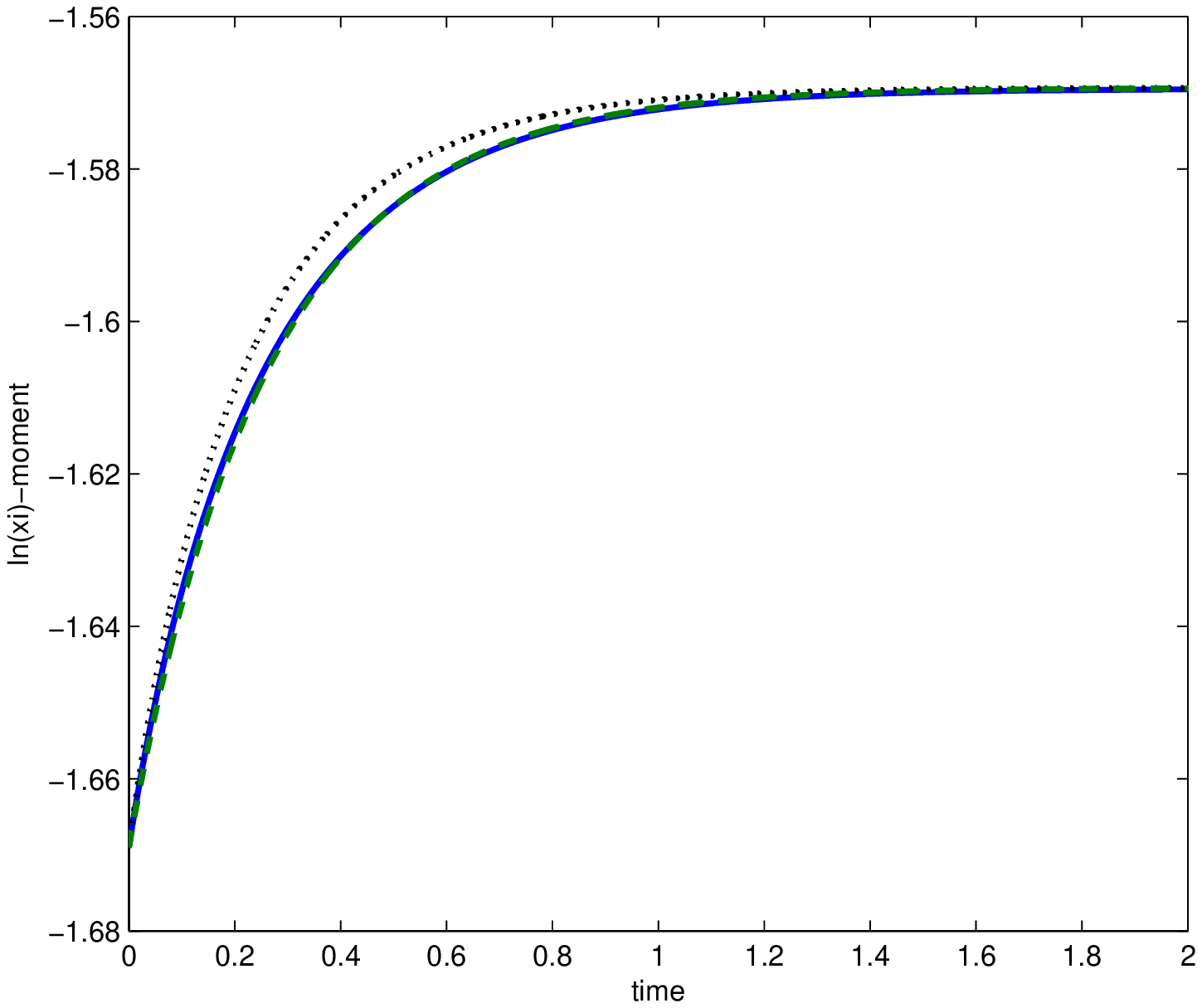}} \\
$\alpha = 2.5$ & $\alpha = 3.0$
\end{tabular}\par}
\caption{Comparison of approximations of the moment $\mmnt{\ln(\xi)}$
obtained by the original \eqref{DNEs}, green dashed curves, and modified \eqref{DNEsB}, red dotted curves, DynMaxEnt methods.
The moment calculated from the solution of the Fokker-Planck equation is plotted with solid blue curves.
Note the different vertical and horizontal axes scales for various values of $\alpha$.
The corresponding relative approximation errors are given in Table \ref{table1}.}
\label{fig:fig1}
\end{figure}

In Fig. \ref{fig:fig1} we plot the time evolution of the $\mmnt{\ln(\xi)}$-moment
of the Fokker-Planck solution $u(t)$ and its approximation obtained
by the DynMaxEnt methods \eqref{DNEs}, \eqref{DNEsB}
for the parameter values $\alpha\in\{1.1, 1.5, 2.5, 3\}$.
Note that since $\alpha > 1$, both the methods \eqref{DNEs}, \eqref{DNEsB} are
applicable. However, we observe that the modified method \eqref{DNEsB}
gives better approximation results. To quantify the approximation error,
we calculate the indicator
\(  \label{error}
   e := \frac{\int_0^T \left( \mmnt{\ln(\xi)}_{u(t)} - \mmnt{\ln(\xi)}_{u_{\alpha^*(t)}} \right)^2 \d t}{\int_0^T \mmnt{\ln(\xi)}_{u(t)}^2 \d t},
\)
where $\mmnt{\ln(\xi)}_{u_{\alpha^*(t)}}$ is the moment calculated by one of the DynMaxEnt
methods  \eqref{DNEs}, \eqref{DNEsB}.
The results for the values $\alpha\in\{1.1, 1.5, 2.5, 3\}$ given in Table \ref{table1}
indeed suggest that the modified method \eqref{DNEsB} provides better approximation
of the moment $\mmnt{\ln(\xi)}_{u(t)}$.
Moreover, we observe that with increasing value of $\alpha$
the approximation properties of both methods seem to improve.

\begin{table}
{\centering \begin{tabular}[h]{|c|cccc|}
\hline
& $\alpha=1.1$ & $\alpha=1.5$ & $\alpha=2.5$ & $\alpha=3.0$ \\
\hline
method \eqref{DNEs} & $1.45\times 10^{-2}$ & $4.09\times 10^{-3}$ & $1.42\times 10^{-3}$ & $1.37\times 10^{-3}$ \\
method \eqref{DNEsB} & $3.78\times 10^{-3}$ & $1.30\times 10^{-3}$ & $3.65\times 10^{-4}$ & $3.41\times 10^{-4}$ \\
\hline
\end{tabular}\par}
\vspace{3mm}

\caption{Relative errors of approximation \eqref{error} of the moment $\mmnt{\ln(\xi)}_{u(t)}$ by the original DynMaxEnt method
\eqref{DNEs}, first row, and its modified version \eqref{DNEsB}, second row. The corresponding plots are given in Fig. \ref{fig:fig1}.}
\label{table1}
\end{table}

\begin{figure}
{\centering \begin{tabular}[h]{cc}
\resizebox*{0.49\linewidth}{!}{\includegraphics{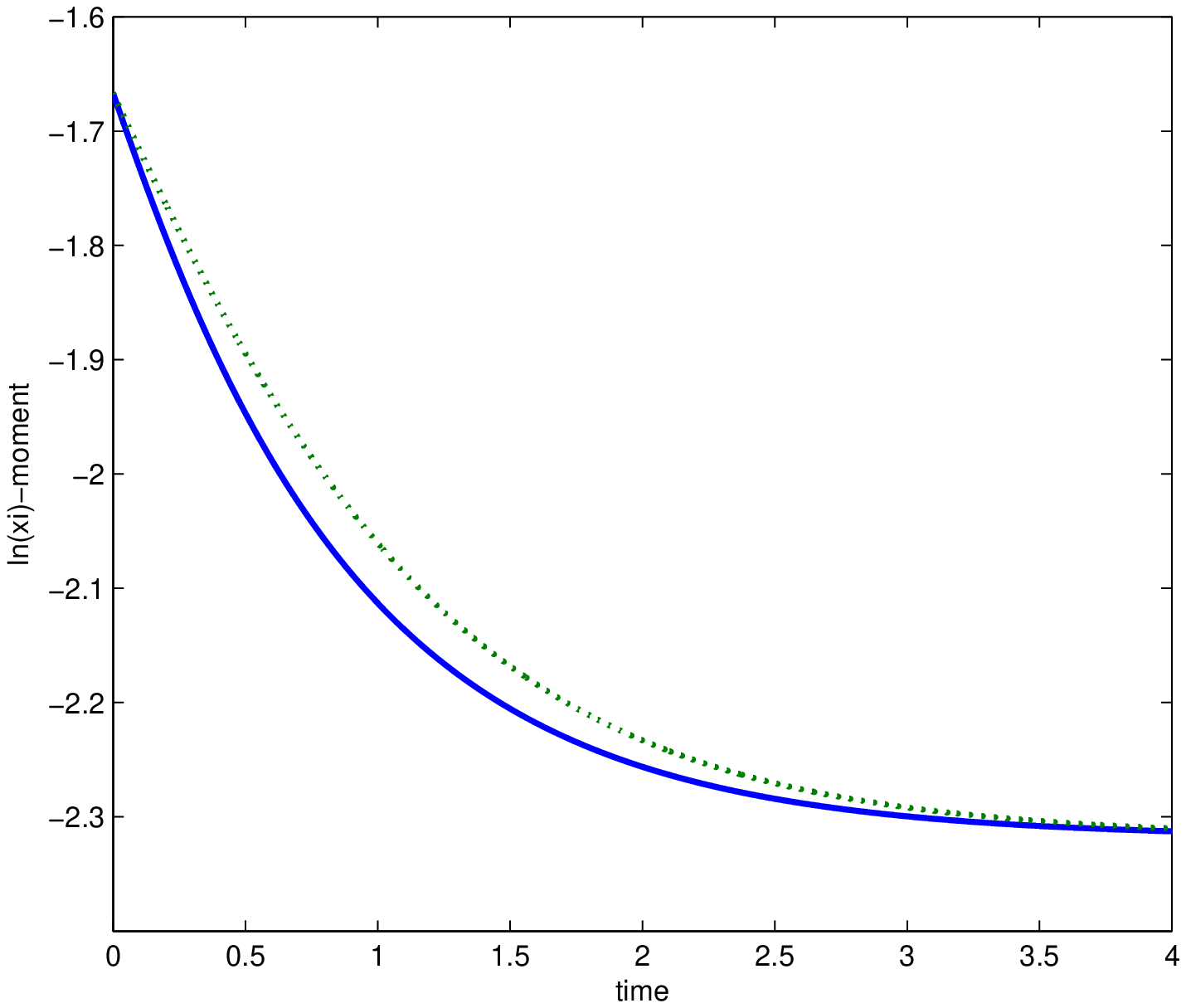}} &
\resizebox*{0.49\linewidth}{!}{\includegraphics{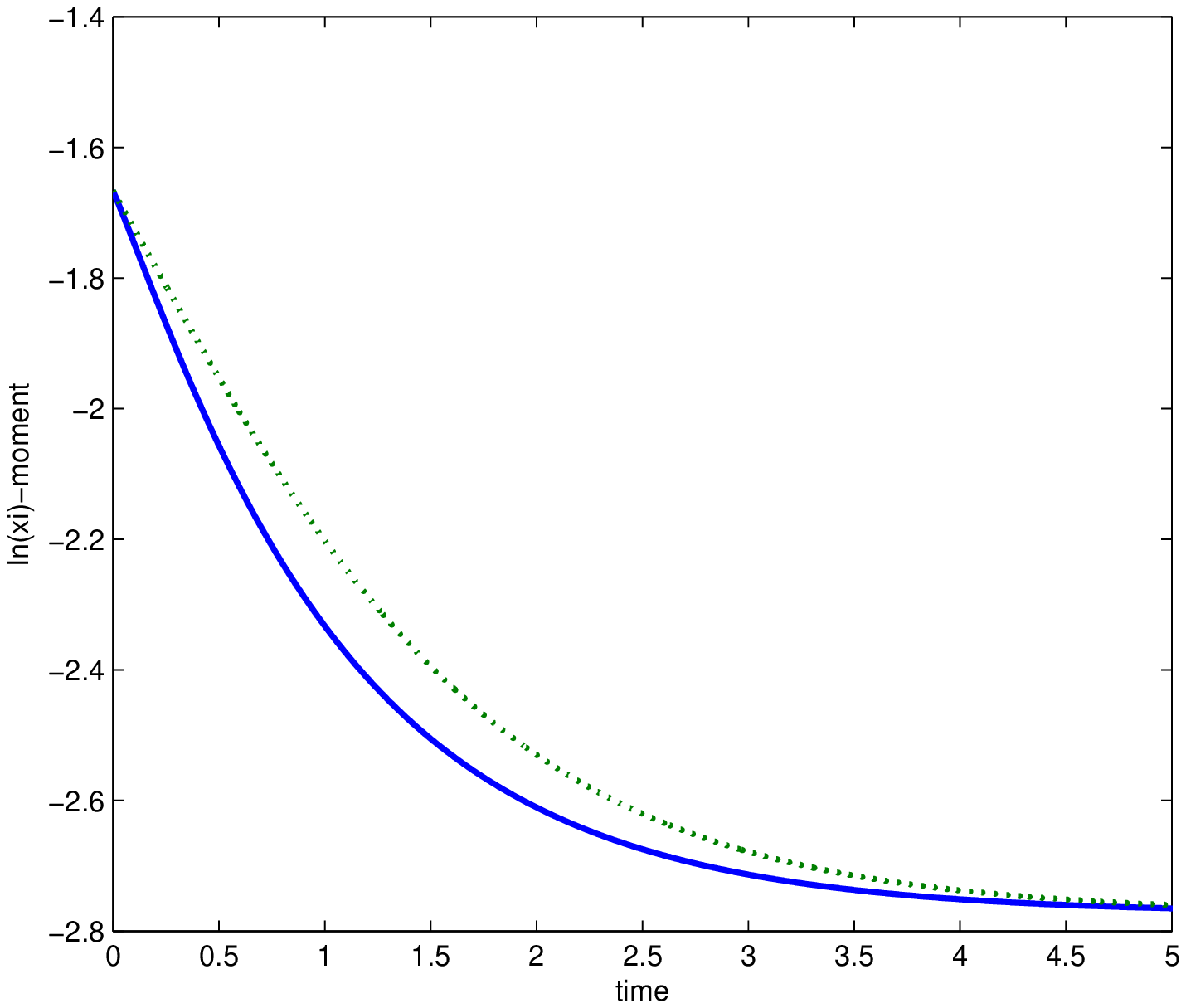}} \\
$\alpha = 0.7$ & $\alpha = 0.5$ \\
\resizebox*{0.49\linewidth}{!}{\includegraphics{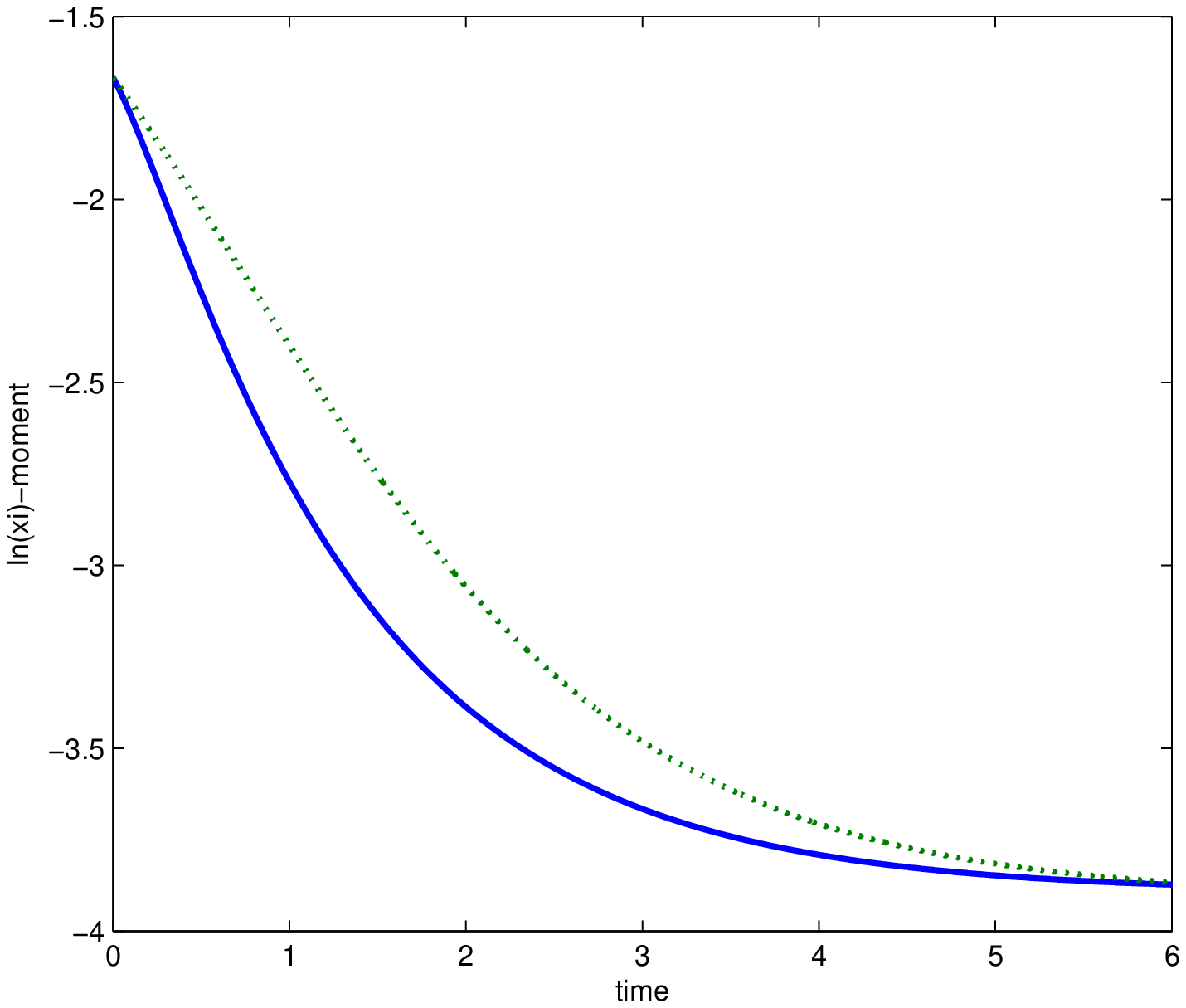}} &
\resizebox*{0.49\linewidth}{!}{\includegraphics{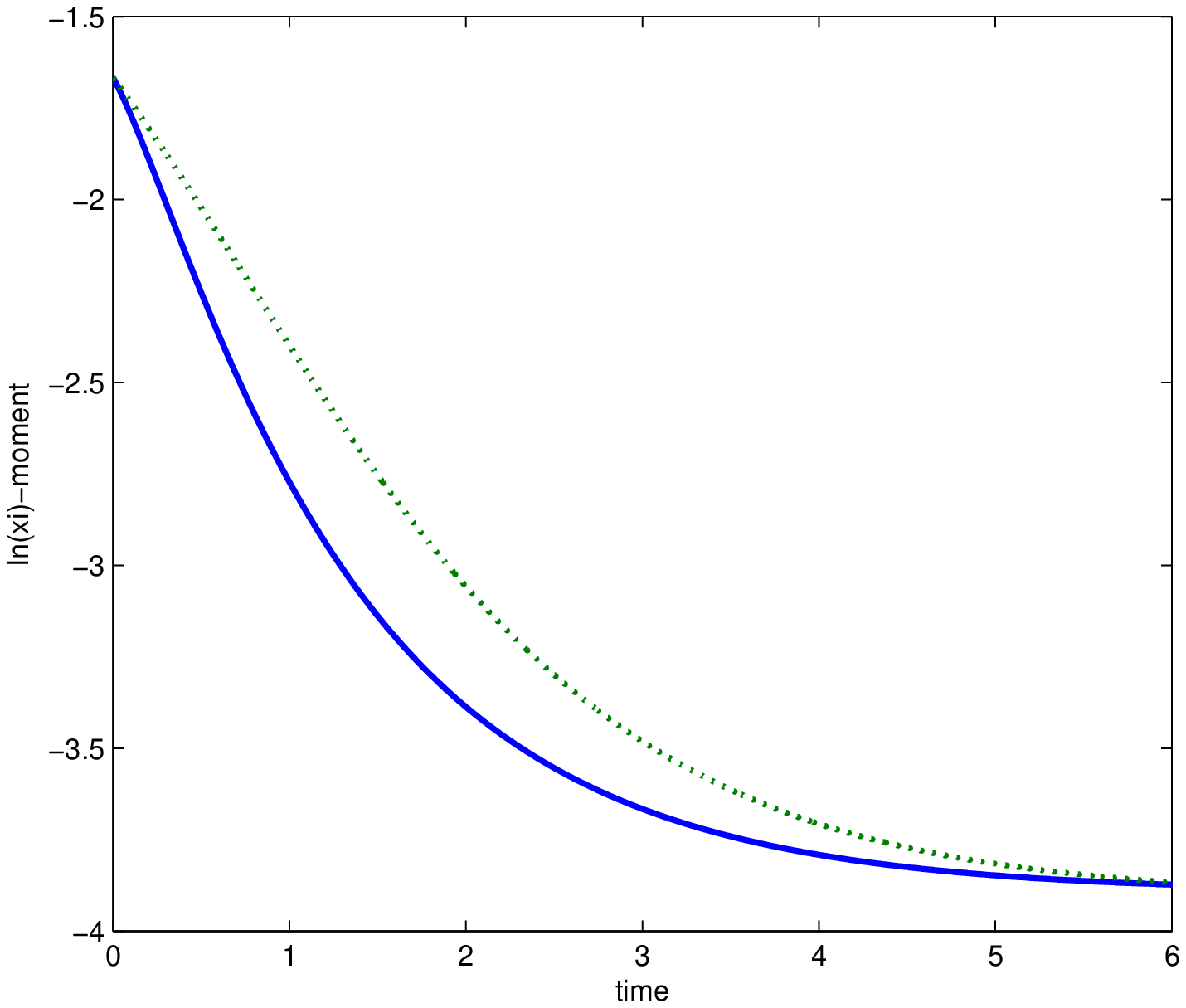}} \\
$\alpha = 0.3$ & $\alpha = 0.2$
\end{tabular}\par}
\caption{Approximations of the moment $\mmnt{\ln(\xi)}$
obtained by the modified \eqref{DNEsB} DynMaxEnt method, red dotted curves, for values of $\alpha < 1$.
The moment calculated from the solution of the Fokker-Planck equation is plotted with solid blue curves.
Note the different vertical and horizontal axes scales for various values of $\alpha$.
The corresponding relative approximation errors are given in Table \ref{table2}.}
\label{fig:fig2}
\end{figure}

In Fig. \ref{fig:fig2} we plot the time evolution of the $\mmnt{\ln(\xi)}$-moment
of the Fokker-Planck solution $u(t)$ and its approximation obtained
by the modified DynMaxEnt method \eqref{DNEsB}
for the parameter values $\alpha\in\{0.7, 0.5, 0.3, 0.2\}$.
Note that the original method \eqref{DNEs} is no longer applicable
since for $\alpha^* < 1$ the term $\mmnt{\xi(\partial_x \ln(\xi))^2}_{u_{\alpha^*}}$
is not finite. Again, we calculate the approximation error \eqref{error}
for the above mentioned valued of $\alpha$ in Table \ref{table2}.
We observe that the approximation worsens for smaller values of $\alpha$.
This is presumably a consequence of the singularity of $u_\alpha$ at $x\in\{0,1\}$
becoming stronger when $\alpha$ approaches zero.
In fact, numerical solution of the Fokker-Planck equation \eqref{eq1} also becomes
more difficult for small values of $\alpha$. For $\alpha<0.2$ our discrete scheme ceases
to provide reliable results. That is why $\alpha=0.2$ is the smallest value
that we take into account.

\begin{table}
{\centering \begin{tabular}[h]{|c|cccc|}
%approximation error \eqref{error}
\hline
& $\alpha=0.7$ & $\alpha=0.5$ & $\alpha=0.3$ & $\alpha=0.2$ \\
\hline
method \eqref{DNEsB} & $1.42\times 10^{-2}$ & $2.81\times 10^{-2}$ & $5.79\times 10^{-2}$ & $8.88\times 10^{-2}$\\
\hline
\end{tabular}\par}
\vspace{3mm}

\caption{Relative errors of approximation \eqref{error} of the moment $\mmnt{\ln(\xi)}_{u(t)}$ by the modified DynMaxEnt method
\eqref{DNEsB}. The corresponding plots are given in Fig. \ref{fig:fig2}.}
\label{table2}
\end{table}

\subsection{Vector case}\label{subsec:num2}
Finally, we consider the more general case
with the function $\bA=\bA(\bx)$ being vector-valued,
$\bA: \Omega_\bx \to \R^k$ with some $k\in\N$.
It has been observed in \cite{dV-Barton12, BTB} that using more moments
(i.e., higher $k$) in general improves the approximation
properties of the DynMaxEnt method.
Inspired by the success of the modified DynMaxEnt method \eqref{DNEsB} demonstrated
in Section \ref{subsec:num1}, we consider an analogous approach
also in the vector case. For this, we employ the idea of deriving a modified DynMaxEnt method
as in Section \ref{subsec:simple}: We multiply the Fokker-Planck equation \eqref{eq1}
by a vector-valued function $\bB: \Omega_\bx \to \R^k$ to be chosen later and integrate by parts,
assuming the boundary terms to vanish. Then, we approximate $u(t)$ by $u_{\balpha^*(t)}$
with the time-dependent vector $\balpha^* = \balpha^*(t)$ and neglect the residual term.
This gives
\(  \label{DynMaxEntB}
    \left( \mmnt{\bB\otimes\bA}_{u_{\balpha^*(t)}}
     - \mmnt{\bB}_{u_{\balpha^*(t)}} \otimes \mmnt{\bA}_{u_{\balpha^*(t)}} \right) \tot{\balpha^*(t)}{t}
     \\
    =  \frac12 \mmnt{\xi \grad_\bx\bB : \grad_\bx\bA}_{u_{\balpha^*(t)}} (\balpha - \balpha^*(t)),
    \nonumber
\)
where $\bB\otimes\bA$ is the $d\times d$ matrix with the $(i,k)$-component $B_i A_k$ and
$\grad_\bx \bB : \grad_\bx \bA$ is the $d\times d$ matrix with the $(i,k)$-component
$\sum_{j=1}^d \partial_{x_j} B_i \partial_{x_j} A_k$.
Clearly, uniform invertibility of the matrix $$\left( \mmnt{\bB\otimes\bA}_{u_{\balpha^*(t)}}
 - \mmnt{\bB}_{u_{\balpha^*(t)}} \otimes \mmnt{\bA}_{u_{\balpha^*(t)}} \right)$$
 is necessary for global solvability of the ODE system \eqref{DynMaxEntB}.
% We were not able to provide an analytic result, however,
 This condition is satisfied in our numerical experiments below.

The goal of this Section is to illustrate the performance of the original \eqref{DynMaxEnt}
and modified \eqref{DynMaxEntB} DynMaxEnt methods for the generic choice $\bA=(\xi',\xi,\ln\xi)$.
For simplicity, we shall still stick to the 1D (single locus) setting $x\in(0,1)$.
Choosing again $\beta = h = 1$ in \eqref{alphaA}, we have
\[
   \aA = - \gamma \xi' + 2 \eta \xi + 2\mu \ln \xi,
\]
with $\balpha = (-\gamma,2\eta,2\mu)$.
The parameters $\gamma, \eta\in\R$ represent the effects of loci on the traits
and $\mu>0$ is the mutation rate.
For the modified DynMaxEnt method \eqref{DynMaxEntB} we choose
$\bB=(\xi',\xi,\xi^2)$. Note that this choice prevents the issue of non-finitness
of the moment $\mmnt{\xi \grad_\bx\bA\otimes\grad_\bx\bA}_{u_{\balpha^*(t)}}$
for $4N\mu < 1$.

We carry out two numerical experiments.
In both simulations we set $N=1$ and the initial condition for the Fokker-Planck equation \eqref{eq1}
to be the stationary distribution \eqref{stat0} with the parameters
$4\mu_0 = 2$, $\eta_0=-1$, $\gamma_0=2$.
As before, the Fokker-Planck equation is discretized in space using the Chang-Cooper
scheme \cite{Chang-Cooper} and forward Euler method in time.
The ODE systems \eqref{DynMaxEnt}, \eqref{DynMaxEntB} are discretized in time
using the forward Euler scheme.

\begin{figure}
{\centering \begin{tabular}[h]{cc}
\resizebox*{0.49\linewidth}{!}{\includegraphics{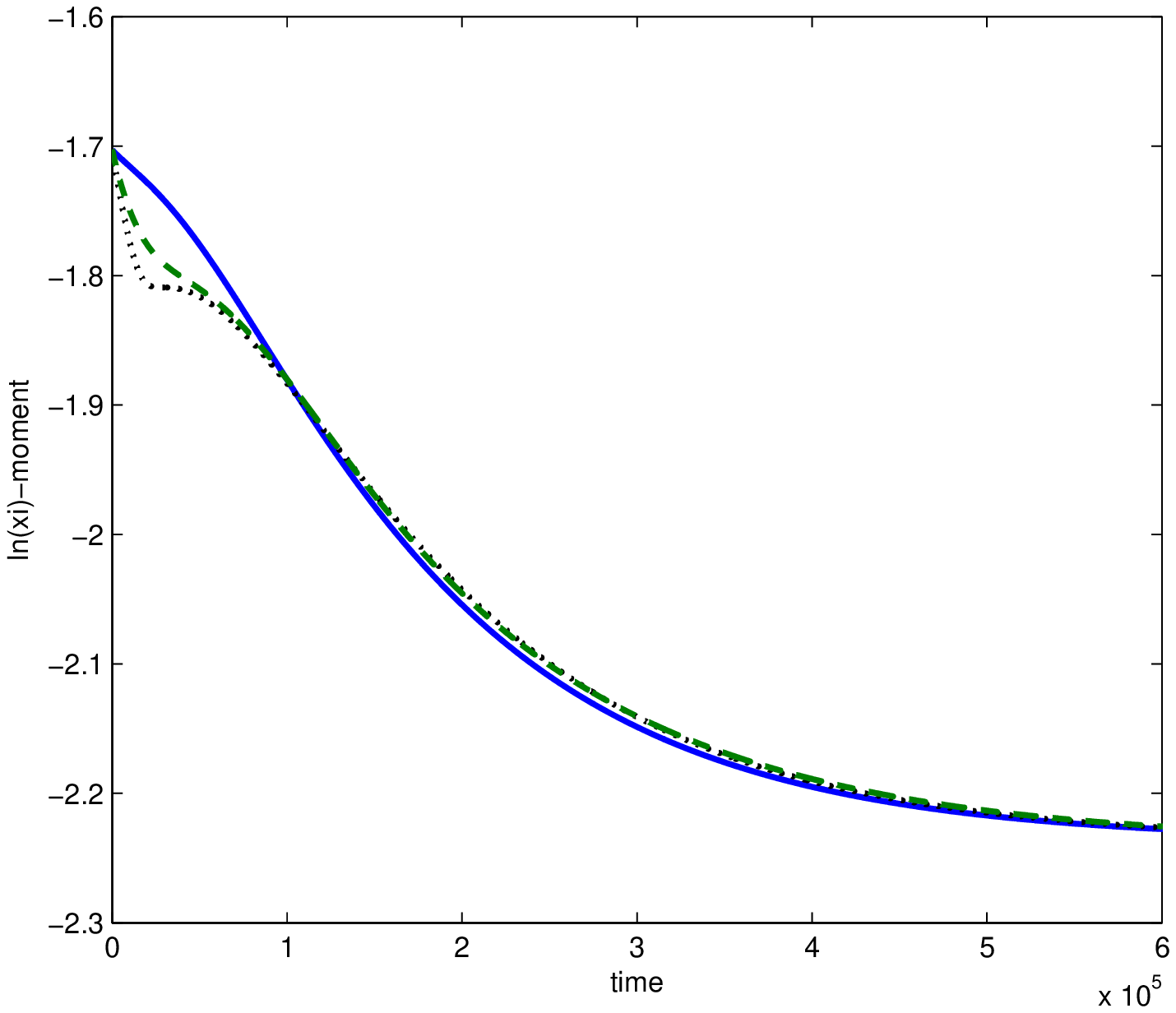}} &
\resizebox*{0.49\linewidth}{!}{\includegraphics{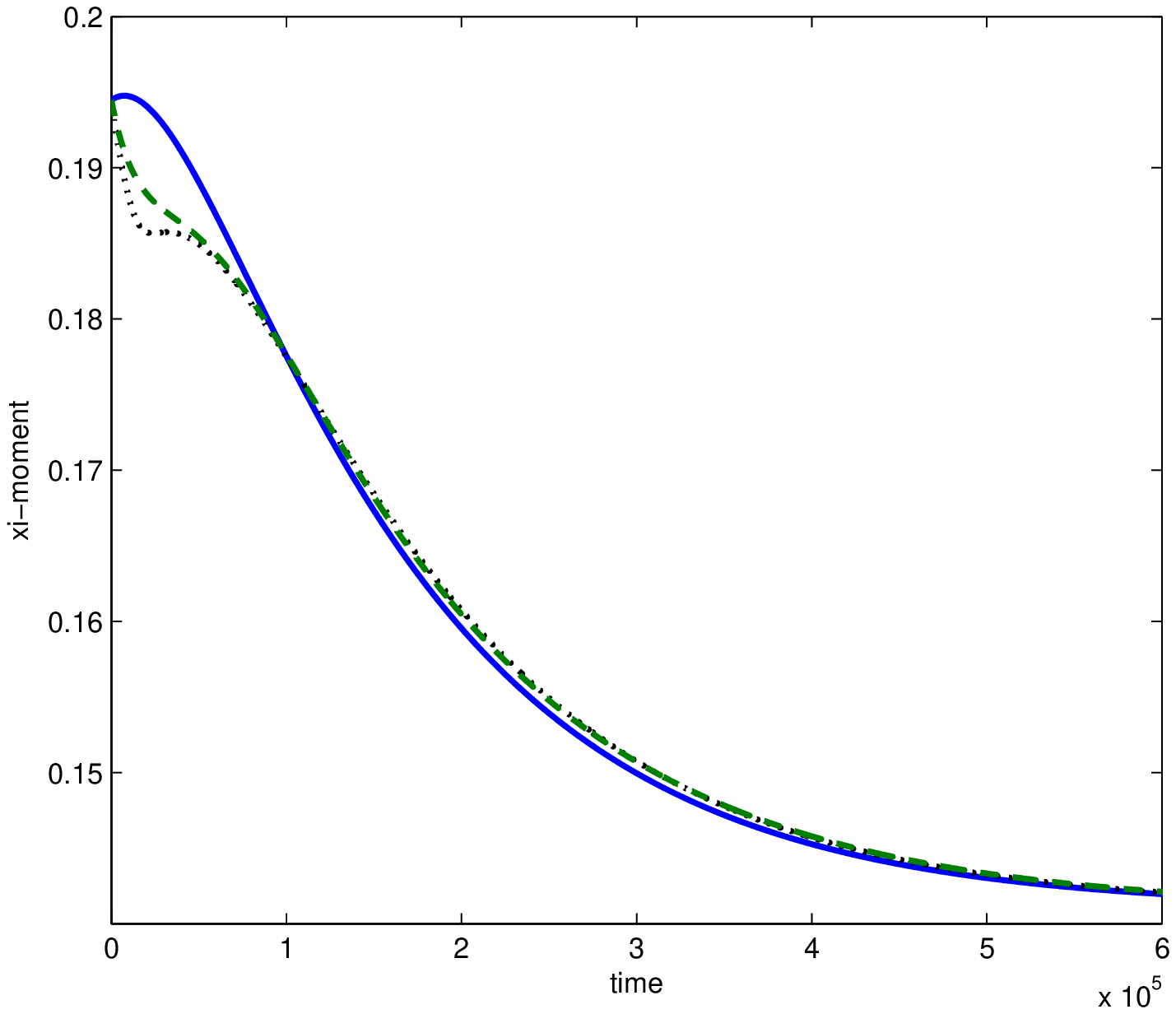}} \\
$\mmnt{\ln(\xi)}$ & $\mmnt{\xi}$ \\
\resizebox*{0.49\linewidth}{!}{\includegraphics{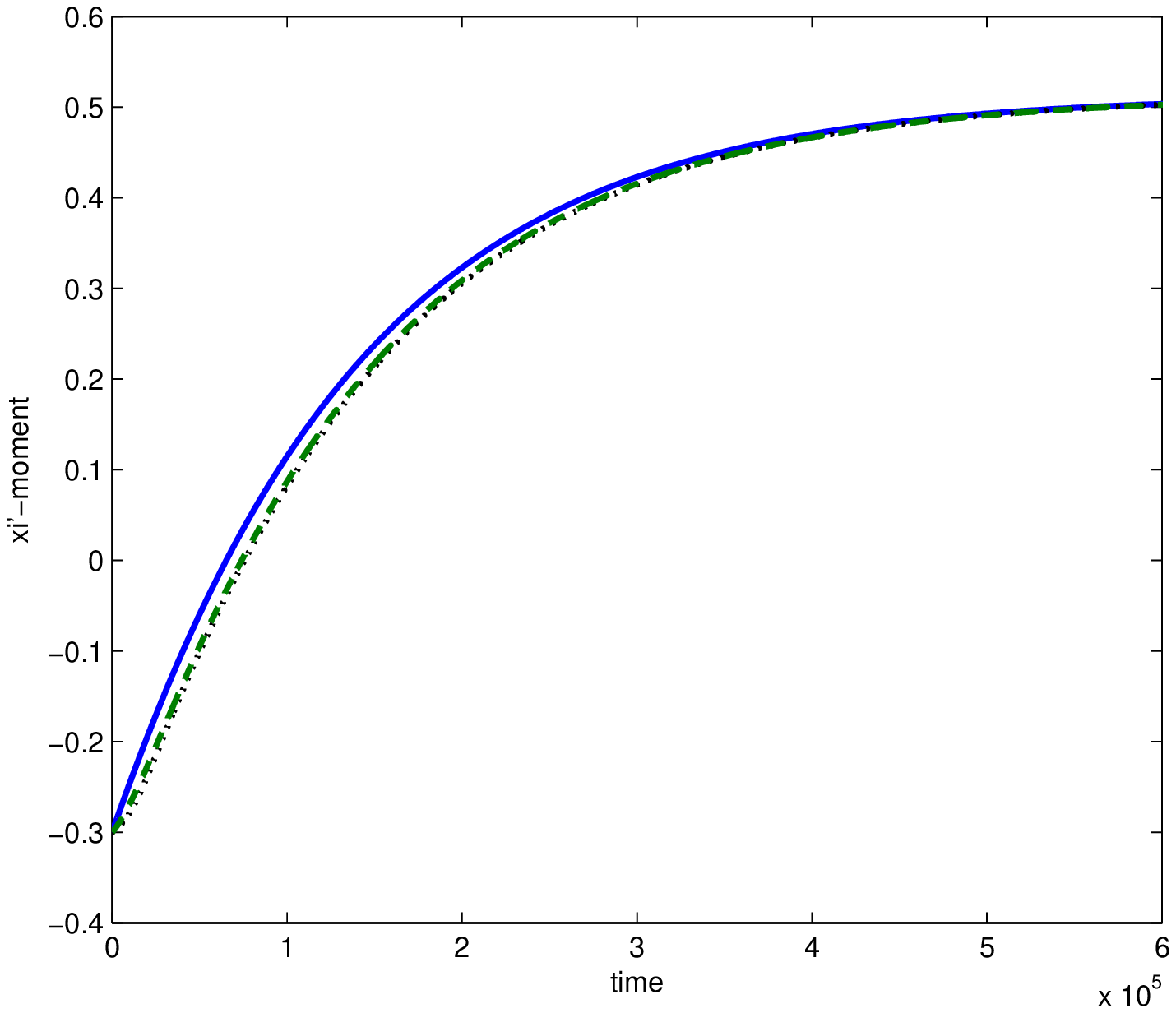}} &
\resizebox*{0.49\linewidth}{!}{\includegraphics{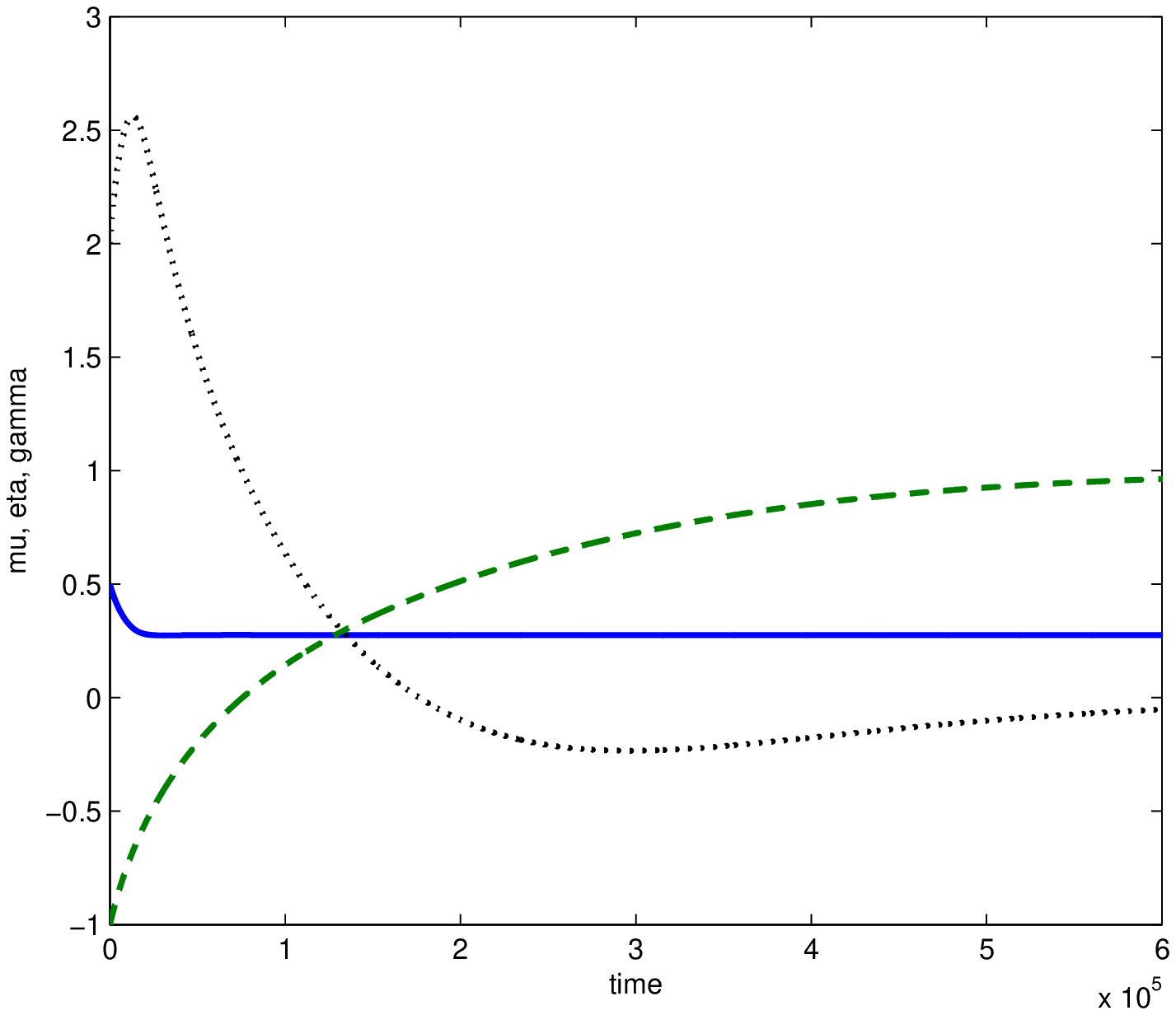}} \\
$\mmnt{\xi'}$ & $\mu$, $\eta$, $\gamma$
\end{tabular}\par}
\caption{Comparison of approximations of the moments $\mmnt{\ln(\xi)}$ (upper left panel),
$\mmnt{\xi}$ (upper right panel) and $\mmnt{\xi'}$ (lower left panel),
obtained by the original \eqref{DynMaxEnt}, black dotted curves, and modified \eqref{DynMaxEntB},
green dashed curves, DynMaxEnt methods.
The moments calculated from the solution of the Fokker-Planck equation \eqref{eq1} is plotted with solid blue curves.
The lower right panel shows the evolution of the parameters $\mu$ (solid blue), $\eta$ (green dashed) and $\gamma$ (black doted)
calculated by the original DynMaxEnt method \eqref{DynMaxEnt}.
The initial and target parameters are $4\mu_0 = 2$, $\eta_0=-1$, $\gamma_0=2$, $4\mu=1.1$, $\eta=1$, $\gamma=0$.
The corresponding relative approximation errors are given in Table \ref{table3}.}
\label{fig:fig3}
\end{figure}

For the first experiment we use the parameter values $4\mu=1.1$, $\eta=1$, $\gamma=0$.
This corresponds to the abrupt change of parameters (evolutionary forces)
\[
   4\mu: 2 \mapsto 1.1,\quad \eta: -1 \mapsto 1,\quad \gamma: 2 \mapsto 0.
\]
In Fig. \ref{fig:fig3} we plot the time evolution of the moments
$\mmnt{\ln(\xi)}$, $\mmnt{\xi}$ and $\mmnt{\xi'}$
of the Fokker-Planck solution $u(t)$ and its approximation obtained
by the original and, resp., modified DynMaxEnt methods \eqref{DynMaxEnt}, resp., \eqref{DynMaxEntB}.
Note that in this case both methods \eqref{DynMaxEnt}, \eqref{DynMaxEntB}
are applicable since the moment $\mmnt{\xi \grad_\bx\bA\otimes\grad_\bx\bA}_{u_{\balpha^*(t)}}$
is finite for all $t\geq 0$.
Calculating the error of approximation \eqref{error} for the three moments,
Table \ref{table3}, we observe that the modified method \eqref{DynMaxEntB}
provides slightly more accurate results.

\begin{table}
{\centering \begin{tabular}[h]{|c|ccc|}
%approximation error \eqref{error}
\hline
&  $\mmnt{\ln(\xi)}$ & $\mmnt{\xi}$ & $\mmnt{\xi'}$ \\
\hline
method \eqref{DynMaxEnt} & $9.24\times 10^{-3}$ & $1.30\times 10^{-2}$ & $5.03\times 10^{-2}$  \\
method \eqref{DynMaxEntB} & $6.79\times 10^{-3}$ & $1.01\times 10^{-2}$ & $4.14\times 10^{-2}$ \\
\hline
\end{tabular}\par}
\vspace{3mm}

\caption{Relative errors of approximation \eqref{error} of the moments $\mmnt{\ln(\xi)}$,
$\mmnt{\xi}$ and $\mmnt{\xi'}$ by the original DynMaxEnt method
\eqref{DynMaxEnt}, first row, and its modified version \eqref{DynMaxEntB}, second row.
The initial and target parameters are $4\mu_0 = 2$, $\eta_0=-1$, $\gamma_0=2$, $4\mu=1.1$, $\eta=1$, $\gamma=0$.
The corresponding plots are given in Fig. \ref{fig:fig3}.}
\label{table3}
\end{table}

\begin{figure}
{\centering \begin{tabular}[h]{cc}
\resizebox*{0.49\linewidth}{!}{\includegraphics{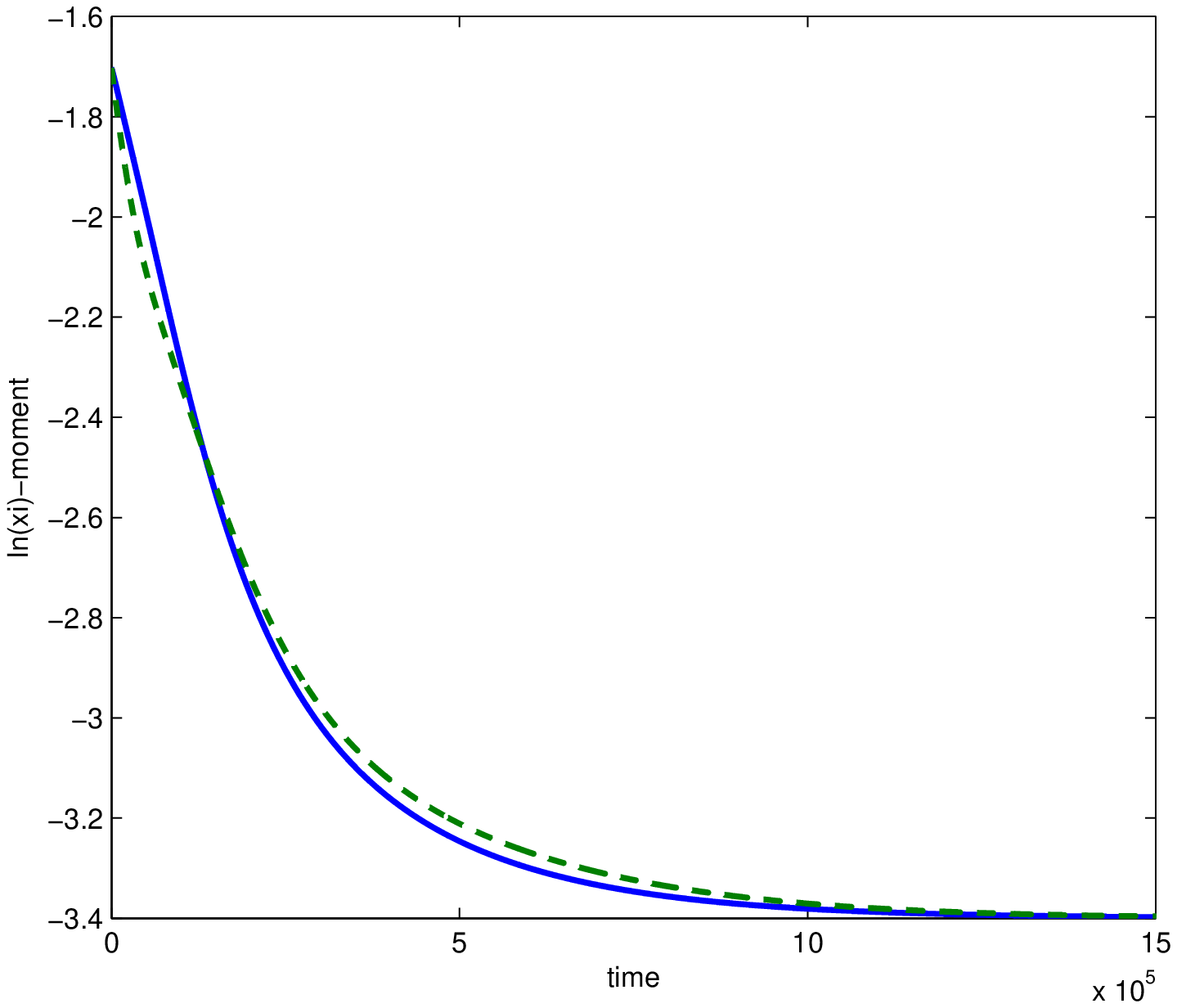}} &
\resizebox*{0.49\linewidth}{!}{\includegraphics{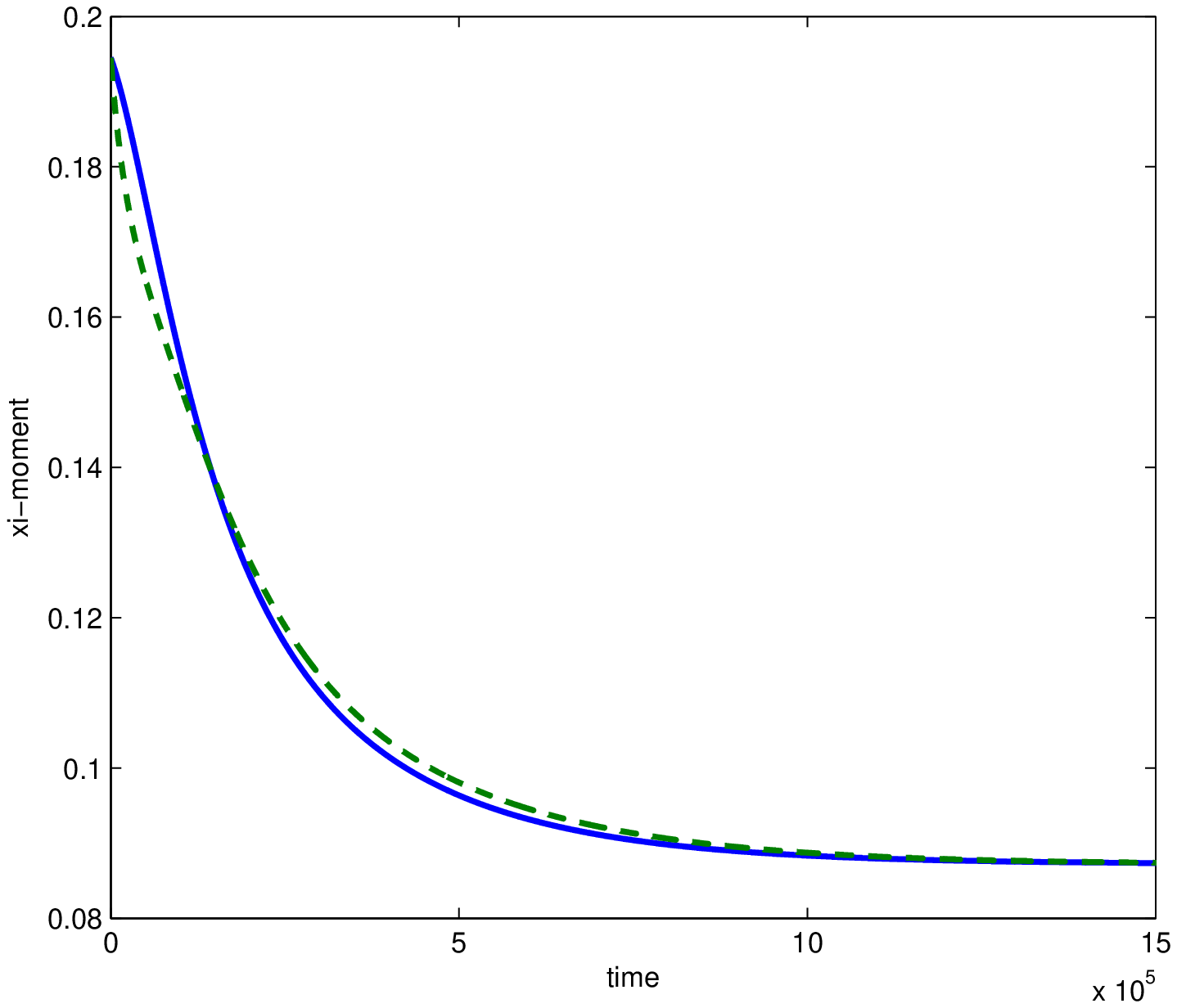}} \\
$\mmnt{\ln(\xi)}$ & $\mmnt{\xi}$ \\
\resizebox*{0.49\linewidth}{!}{\includegraphics{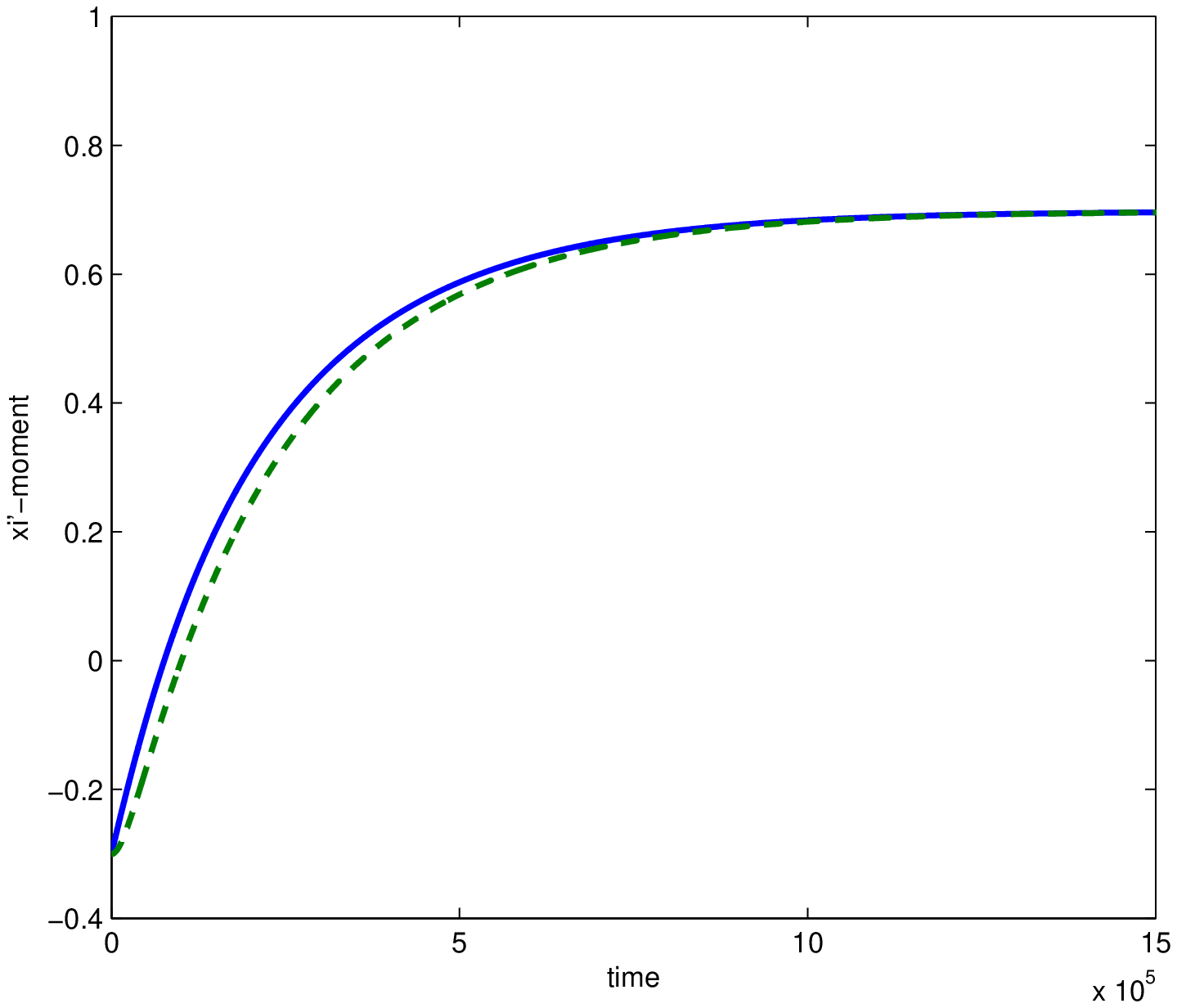}} &
\resizebox*{0.49\linewidth}{!}{\includegraphics{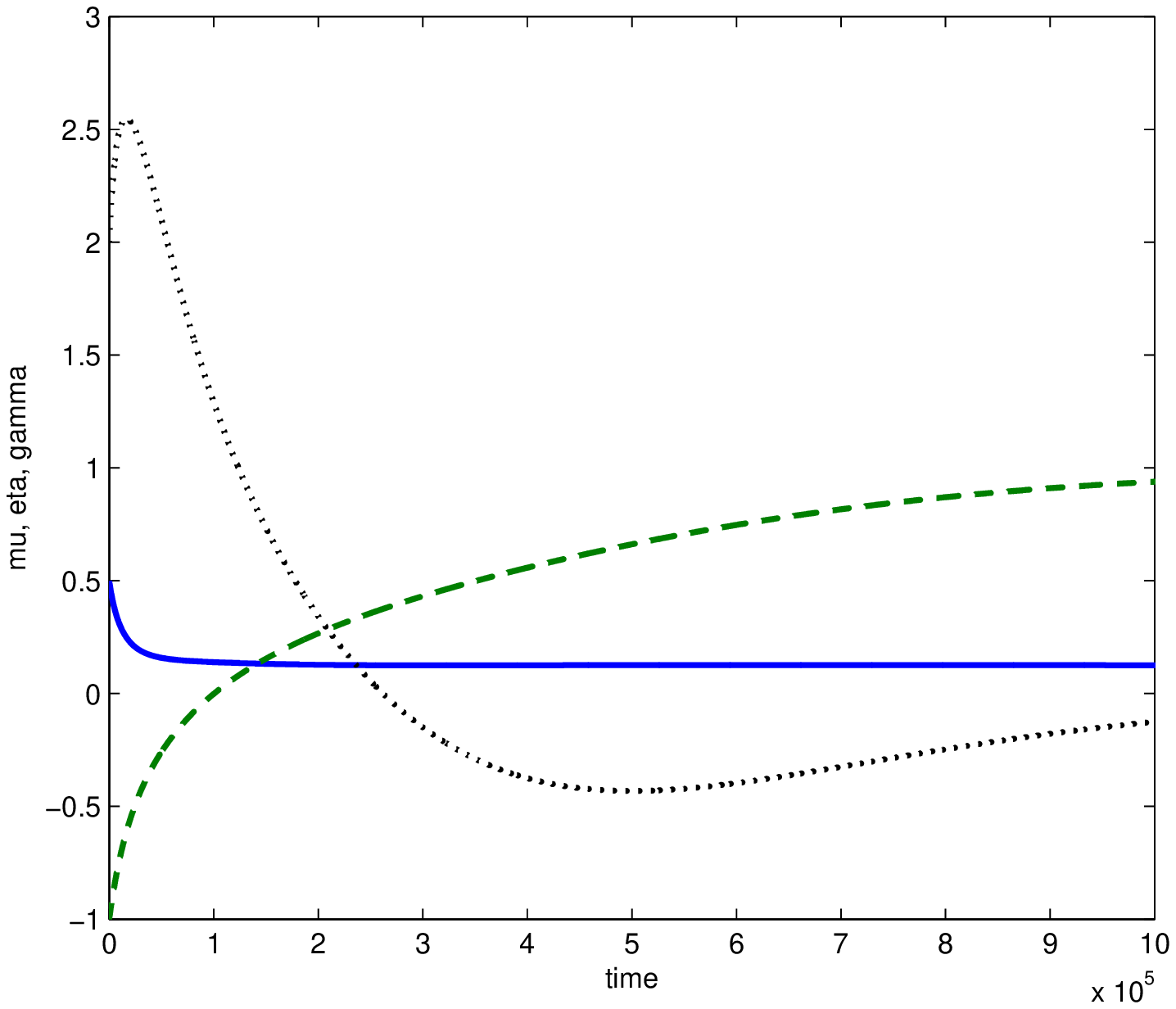}} \\
$\mmnt{\xi'}$ & $\mu$, $\eta$, $\gamma$
\end{tabular}\par}
\caption{Approximations of the moments $\mmnt{\ln(\xi)}$ (upper left panel),
$\mmnt{\xi}$ (upper right panel) and $\mmnt{\xi'}$ (lower left panel),
obtained by the modified \eqref{DynMaxEntB} DynMaxEnt method (green dashed curves).
The moments calculated from the solution of the Fokker-Planck equation \eqref{eq1} is plotted with solid blue curves.
The lower right panel shows the evolution of the parameters $\mu$ (solid blue), $\eta$ (green dashed) and $\gamma$ (black doted).
The initial and target parameters are $4\mu_0 = 2$, $\eta_0=-1$, $\gamma_0=2$, $4\mu=0.5$, $\eta=1$, $\gamma=0$.
The corresponding relative approximation errors are given in Table \ref{table4}.}
\label{fig:fig4}
\end{figure}

For our second experiment we use the parameter values $4\mu=0.5$, $\eta=1$, $\gamma=0$.
This corresponds to the rapid change of evolutionary forces
\[
   4\mu: 2 \mapsto 0.5,\quad \eta: -1 \mapsto 1,\quad \gamma: 2 \mapsto 0.
\]
Note that in this case the original method \eqref{DynMaxEnt} is not applicable any more
since the moment $\mmnt{\xi \grad_\bx\bA\otimes\grad_\bx\bA}_{u_{\balpha^*(t)}}$
is not defined for $4\mu < 1$.
In Fig. \ref{fig:fig4} we plot the time evolution of the moments
$\mmnt{\ln(\xi)}$, $\mmnt{\xi}$ and $\mmnt{\xi'}$
of the Fokker-Planck solution $u(t)$ and its approximation obtained
by the modified DynMaxEnt method \eqref{DynMaxEntB}.
On the other hand, the results presented in Fig. \ref{fig:fig4} and Table \ref{table4} indicate
that the modified DynMaxEnt method \eqref{DynMaxEntB} provides a reasonably
good approximation of the three moments.

\begin{table}
{\centering \begin{tabular}[h]{|c|ccc|}
%approximation error \eqref{error}
\hline
&  $\mmnt{\ln(\xi)}$ & $\mmnt{\xi}$ & $\mmnt{\xi'}$ \\
\hline
method \eqref{DynMaxEntB} & $2.45\times 10^{-2}$ & $2.55\times 10^{-2}$ & $1.24\times 10^{-1}$ \\
\hline
\end{tabular}\par}
\vspace{3mm}

\caption{Relative errors of approximation \eqref{error} of the moments $\mmnt{\ln(\xi)}$,
$\mmnt{\xi}$ and $\mmnt{\xi'}$ by the modified DynMaxEnt method \eqref{DynMaxEntB}.
The initial and target parameters are $4\mu_0 = 2$, $\eta_0=-1$, $\gamma_0=2$, $4\mu=0.5$, $\eta=1$, $\gamma=0$.
The corresponding plots are given in Fig. \ref{fig:fig4}.}
\label{table4}
\end{table}

$\,$\\
\noindent\textbf{Acknowledgments.}
We thank Nicholas Barton (IST Austria) for his useful comments and suggestions.
JH and PM are funded by KAUST baseline funds and grant no. 1000000193.

%%%%%%%%%%%%%%%%%%

%%%%%%%%%%%%%%%%%%


\begin{thebibliography}{10}

\bibitem{ALT}
D. Albanez, H. Nussenzveig Lopes, and E. Titi:
\emph{Continuous data assimilation for the three-dimensional Navier-Stokes--$\alpha$ model.}
Asymptotic Analysis 97 (2016), 139--164.

\bibitem{ATKZ}
M. Altaf, E. Titi, O. Knio, L. Zhao, M. McCabe, and I. Hoteit:
\emph{Downscaling the 2D B\'enard Convection Equations Using Continuous Data Assimilation.}
Computational Geosciences (to appear, 2017).

\bibitem{AMTU} A. Arnold, P. A. Markowich, G. Toscani, and A. Unterreiter:
\emph{On convex Sobolev inequalities and the rate of convergence to equilibrium for Fokker-Planck type equations.}
Comm. PDE 26 (2001), 43--100.

\bibitem{Barton09}
N. Barton, and H. de Vladar: \emph{Statistical mechanics and the evolution of polygenic quantitative traits.}
Genetics 181 (2009), 997--1011.

\bibitem{Bialek}
W. Bialek, A. Cavagna, I. Giardina, T. Mora, E. Silvestri et al.:
\emph{Statistical mechanics for natural flocks of birds.}
Proc. Natl. Acad. Sci. USA 109 (2012), 4786--4791.

\bibitem{BTB}
K. Bo{\md}ov\'a, G. Tkacik, and N. Barton: \emph{A General Approximation for the Dynamics of Quantitative Traits.}
Genetics, Vol. 202 (2016), 1523--1548.

\bibitem{Chalub1}
F. Chalub and M. Souza: \emph{From discrete to continuous evolution models: A unifying approach to
drift-diffusion and replicator dynamics.}
Theoretical Population Biology 76 (2009), 268--277.

\bibitem{Chalub2}
F. Chalub and M. Souza: \emph{A non-standard evolution problem arising in population genetics.}
Comm. Math. Sci. 7 (2009), 489--502.

\bibitem{Chang-Cooper}
J. S. Chang, and G. Cooper: \emph{A practical difference scheme for Fokker-Planck equation.}
Journal of Computational Physics 6 (1970), 1--16.

\bibitem{dV-Barton12}
H. de Vladar and N. H. Barton: \emph{The statistical mechanics of a polygenic character under stabilizing selection, mutation and drift.}
J. R. Soc. Interface 8 (2012), 720--739.

\bibitem{Ewens12}
W. J. Ewens: \emph{Mathematical Population Genetics 1: Theoretical Introduction.}
Interdisciplinary Applied Mathematics, Vol. 27. Springer, New York, 2012.

\bibitem{FLT}
A. Farhat, E. Lunasin, and E. Titi: \emph{Continuous data assimilation algorithm for a 2D B\'enard convection
through horizontal velocity measurements alone.}
Journal of Nonlinear Science (online first, 2017).

\bibitem{GOT}
M. Gesho, E. Olson, and E. Titi: \emph{A Computational Study of a Data Assimilation Algorithm for the
Two-dimensional Navier--Stokes Equations.}
Communications in Computational Physics 19 (2016), 1094--1110.

\bibitem{Hick}
P. Hick, and G. Stevens: \emph{Approximate solutions to the cosmic ray transport equation the maximum entropy method.}
Astron. Astrophys. 172 (1987), 350--358. 

\bibitem{Kimura55a}
M. Kimura: \emph{Solution of a process of random genetic drift with a continuous model.}
Proc. Natl. Acad. Sci. USA 41 (1955), 144.

\bibitem{Kimura55b}
M. Kimura: \emph{Stochastic processes and distribution of gene frequencies under natural selection.}
Cold Spring Harb. Symp. Quant. Biol. 20 (1955), 33--53.

\bibitem{MTT}
P. Markowich, E. Titi, and S. Trabelsi:
\emph{Continuous data assimilation for the three-dimensional Brinkman-Forchheimer-extended Darcy Model.}
Nonlinearity 29 (2016).

\bibitem{MT}
C. Mondaini and E.S. Titi: \emph{Postprocessing Galerkin method applied to a data assimilation algorithm:
a uniform in time error estimate.}
SIAM Journal on Numerical Analysis (submitted, 2017).

\bibitem{Mora}
T. Mora, A. M. Walczak, W. Bialek, and C. G. Callan:
\emph{Maximum entropy models for antibody diversity.}
Proc. Natl. Acad. Sci. USA 107 (2010), 5405--5410.

\bibitem{Pazy}
A. Pazy: \emph{Semigroups of Linear Operators and Applications to Partial Differential Equations.}
Applied Mathematical Sciences 44, Springer-Verlag, Berlin, 1983.

\bibitem{Plastino}
A.R. Plastino, H. Miller, and A. Plastino: \emph{Minimum Kullback entropy approach to the Fokker--Planck equation.}
Phys. Rev. E 56 (1997), 3927--3934.

\bibitem{Shapiro97}
A. Pr\"ugel-Bennett, and J. Shapiro: \emph{An analysis of genetic algorithms using statistical mechanics.}
Physica D 104 (1997), 75--114. 

\bibitem{Shapiro01}
M. Rattray, and J. L. Shapiro: \emph{Cumulant dynamics of a population under multiplicative selection, mutation, and drift.}
Theor. Popul. Biol. 60 (2001), 17--31.

\bibitem{Schneidman}
E. Schneidman, M. J. Berry, R. Segev, and W. Bialek:
\emph{Weak pairwise correlations imply strongly correlated network states in a neural population.}
Nature 440 (2006), 1007--1012.

\bibitem{Simon}
B. Simon: \emph{Schr\"odinger operators with purely discrete spectrum.}
Methods of Functional Analysis and Topology 15 (2009), 61--66.

\bibitem{Tkacik}
G. Tka\v{c}ik, O. Marre, D. Amodei, E. Schneidman, W. Bialek et al.:
\emph{Searching for collective behavior in a large network of sensory neurons.}
PLoS Comut. Biol. 10 (2014), e1003408.

%\bibitem{Wakeley08}
%J. Wakeley: \emph{Coalescent Theory: An Introduction.}
%Roberts \& Company, Greenwood Village, CO, 2008.

\bibitem{Weigt}
M. Weigt, R. A. White, H. Szurmant, J. A. Hoch, and T. Hwa:
\emph{Identification of direct residue contacts in protein-protein interaction by message passing.}
Proc. Natl. Acad. Sci. USA 106 (2009), 67--72.

\end{thebibliography}
\end{document}